\documentclass[11pt,english]{article}
\usepackage[T1]{fontenc}
\usepackage[latin9]{inputenc}
\usepackage{geometry}
\usepackage{color}
\usepackage{units}
\usepackage{url}
\usepackage{amsmath}
\usepackage{amssymb}

\makeatletter
\usepackage{babel}

\makeatother

\usepackage{babel}
\begin{document}

\title{On the spectra of hypermatrix direct sum\\
 and Kronecker products constructions.}

\author{Yuval Filmus \thanks{Computer Science Department, Technion - Israel Institute of Technology},
Edinah K. Gnang \thanks{Department of Mathematics, Purdue University, email: egnang@purdue.edu}}
\maketitle
\begin{abstract}
We extend to hypermatrices definitions and theorem from matrix theory.
Our main result is an elementary derivation of the spectral decomposition
of hypermatrices generated by arbitrary combinations of Kronecker
products and direct sums of cubic side length $2$ hypermatrices.
The method is based on a generalization of Parseval's identity. We
use this general formulation of Parseval's identity to introduce hypermatrix
Fourier transforms and discrete Fourier hypermatrices. We extend to
hypermatrices a variant of the Gram\textendash Schmidt orthogonalization
process as well as Sylvester's classical Hadamard matrix construction.
We conclude the paper with illustrations of spectral decompositions
of adjacency hypermatrices of finite groups and a short proof of the
hypermatrix formulation of the Rayleigh quotient inequality.
\end{abstract}

\section{Introduction}

\emph{Hypermatrices} are multidimensional arrays of complex numbers
which generalize matrices. Formally, we define a hypermatrix to be
a finite set of complex numbers indexed by distinct elements of some
fixed integer Cartesian product set of the form 
\[
\left\{ 0,1,2,\cdots,n_{1}\right\} \times\left\{ 0,1,2,\cdots,n_{2}\right\} \times\cdots\times\left\{ 0,1,2,\cdots,n_{m}\right\} .
\]
Such a hypermatrix is said to be of order $m$ and of size $\left(n_{1}+1\right)\times\left(n_{2}+1\right)\times\cdots\times\left(n_{m}+1\right)$.
The hypermatrix is said to be cubic and of side length $\left(n+1\right)$
if $n_{1}=n_{2}=\cdots=n_{m}=n$. In particular, matrices are second
order hypermatrices. Hypermatrix algebras arise from attempts to extend
to hypermatrices classical matrix algebra concepts and algorithms
\cite{MB94,GKZ,RK,GER}. Hypermatrices are common occurrences in applications
relating to computer science, statistics and physics. In these applications
hypermatrices are often embedded into multilinear forms associated
with objective functions to be minimized or maximized. While many
hypermatrix algebras have been proposed in the tensor/hypermatrix
literature \cite{Lim2013}, our discussion here focuses on the Bhattacharya-Mesner
(BM) algebra first developed in \cite{MB90,MB94} and the general
BM algebra first proposed in \cite{GER}. The general BM product is
of interest because it encompasses as special cases many other hypermatrix
products discussed in the literature including the usual matrix product,
the Segre outer product, the contraction product, the higher order
singular value decomposition, and the multilinear matrix multiplication. 

The study of structured dense matrices such as discrete Fourier matrices,
Hadamard matrices, Hankel matrices, Hessian matrices, Vandermonde
matrices, Wronksian matrices, as well as matrix direct sum and Kronecker
product constructions play important roles in the applications of
linear algebra methods to other disciplines. Many established results
concerning such matrices draw heavily on matrix spectral analysis
toolkits which include techniques derived from matrix spectral decompositions,
Fourier-Hadamard\textendash Rademacher\textendash Walsh transforms,
the Parseval's identity, the Gram\textendash Schmidt orthogonalization
process and the Rayleigh quotient inequality. The development of hypermatrix
spectral analysis toolkits constitutes the main obstacle to extending
these results to hypermatrices. The present work aims to add to hypermatrix
spectral analysis toolkits. Our main result is a constructive method
for obtaining spectral decomposition of hypermatrix direct sum and
Kronecker product constructions. The method is based in part on a
generalization of Parseval's identity. We use this general formulation
of Parseval's identity to introduce hypermatrix Fourier transforms
and discrete Fourier hypermatrices. We extend to hypermatrices a variant
of the Gram\textendash Schmidt orthogonalization process as well as
Sylvester's classical Hadamard matrix construction. We conclude the
paper with illustrations of spectral decompositions of adjacency hypermatrices
of finite groups and a short proof of the hypermatrix Rayleigh quotient
inequality.

This article is accompanied by an extensive and actively maintained
Sage \cite{sage} symbolic hypermatrix algebra package which implements
the various features of the general BM algebra. The package is made
available at the link \url{https://github.com/gnang/HypermatrixAlgebraPackage}\\
\\
\emph{Acknowledgement}: We would like to thank Andrei Gabrielov for
providing guidance while preparing this manuscript. We would like
to thank Vladimir Retakh and Ahmed Elgammal for patiently introducing
us to the theory of hypermatrices. We are grateful to Doron Zeilberger
and Saugata Basu for helpful suggestions. We would like to thank Ori
Parzanchevski for helping get the hypermatrix package started. The
second author was supported by the National Science Foundation, and
is grateful for the hospitality of the Institute for Advanced Study
and the Department of Mathematics at Purdue University.

\section{Overview of the BM algebra}

The Bhattacharya-Mesner product, or BM product for short, was first
developed in \cite{MB90,MB94}. The BM product provides a natural
generalization to the matrix product. The BM product of second order
hypermatrices corresponds to the usual matrix product. For notational
consistency, we will on occasion use the notation $\mbox{Prod}\left(\mathbf{A}^{(1)},\mathbf{A}^{(2)}\right)$
to refer to the matrix product $\mathbf{A}^{(1)}\cdot\mathbf{A}^{(2)}$.
The BM product is best introduced to the unfamiliar reader by first
describing the BM product of third and fourth order hypermatrices.
Note that the BM product of second order hypermatrices is a binary
operation, the BM product of third order hypermatrices is a ternary
operation, the BM product of fourth order hypermatrices takes four
operands, and so on.

The BM product of third order hypermatrices $\mathbf{A}^{(1)}$, $\mathbf{A}^{(2)}$
and $\mathbf{A}^{(3)}$ , denoted $\mbox{Prod}\left(\mathbf{A}^{(1)},\mathbf{A}^{(2)},\mathbf{A}^{(3)}\right)$,
is defined if 
\[
\mathbf{A}^{(1)}\mbox{ is }n_{1}\times\textcolor{red}{k}\times n_{3},\:\mathbf{A}^{(2)}\mbox{ is }n_{1}\times n_{2}\times\textcolor{red}{k}\,\mbox{ and }\mathbf{A}^{(3)}\mbox{ is }\textcolor{red}{k}\times n_{2}\times n_{3}.
\]
The result $\mbox{Prod}\left(\mathbf{A}^{(1)},\mathbf{A}^{(2)},\mathbf{A}^{(3)}\right)$
will be of size $n_{1}\times n_{2}\times n_{3}$, and specified entry-wise
by 
\[
\left[\mbox{Prod}\left(\mathbf{A}^{(1)},\mathbf{A}^{(2)},\mathbf{A}^{(3)}\right)\right]_{i_{1},i_{2},i_{3}}=\sum_{0\le\textcolor{red}{j}<k}a_{i_{1}\,\textcolor{red}{j}\,i_{3}}^{(1)}\,a_{i_{1}\,i_{2}\,\textcolor{red}{j}}^{(2)}\,a_{\textcolor{red}{j}\,i_{2}\,i_{3}}^{(3)}.
\]
Similarly, the BM product of fourth order hypermatrices $\mathbf{A}^{(1)}$
,$\mathbf{A}^{(2)}$ ,$\mathbf{A}^{(3)}$ and $\mathbf{A}^{(4)}$,
denoted $\mbox{Prod}\left(\mathbf{A}^{(1)},\mathbf{A}^{(2)},\mathbf{A}^{(3)},\mathbf{A}^{(4)}\right)$
is defined if
\[
\mathbf{A}^{(1)}\mbox{ is }n_{1}\times\textcolor{red}{k}\times n_{3}\times n_{4},\,\mathbf{A}^{(2)}\mbox{ is }n_{1}\times n_{2}\times\textcolor{red}{k}\times n_{4},
\]
\[
\mathbf{A}^{(3)}\mbox{ is }n_{1}\times n_{2}\times n_{3}\times\textcolor{red}{k},\mbox{ and }\mathbf{A}^{(4)}\mbox{ is }\textcolor{red}{k}\times n_{2}\times n_{3}\times n_{4}.
\]
The result $\mbox{Prod}\left(\mathbf{A}^{(1)},\mathbf{A}^{(2)},\mathbf{A}^{(3)},\mathbf{A}^{(4)}\right)$
will be of size $n_{1}\times n_{2}\times n_{3}\times n_{4}$ and specified
entry-wise by 
\[
\left[\mbox{Prod}\left(\mathbf{A}^{(1)},\mathbf{A}^{(2)},\mathbf{A}^{(3)},\mathbf{A}^{(4)}\right)\right]_{i_{1},i_{2},i_{3},i_{4}}=\sum_{0\le\textcolor{red}{j}<k}a_{i_{1}\,\textcolor{red}{j}\,i_{3}i_{4}}^{(1)}\,a_{i_{1}i_{2}\,\textcolor{red}{j}\,i_{4}}^{(2)}\,a_{i_{1}i_{2}i_{3}\,\textcolor{red}{j}}^{(3)}\:a_{\textcolor{red}{j}\,i_{2}i_{3}i_{4}}^{(4)}.
\]
The reader undoubtedly has already discerned the general pattern,
but for the sake of completeness we express the entries of the BM
product of order $m$ hypermatrices 
\begin{equation}
\left[\mbox{Prod}\left(\mathbf{A}^{(1)},\cdots,\mathbf{A}^{(t)},\cdots,\mathbf{A}^{(m)}\right)\right]_{i_{1},\cdots,i_{t},\cdots,i_{m}}=\sum_{0\le\textcolor{red}{j}<k}a_{i_{1}\,\textcolor{red}{j}\,i_{3}\cdots i_{m}}^{(1)}\cdots a_{i_{1}\cdots i_{t}\,\textcolor{red}{j}\,i_{t+2}\cdots i_{m}}^{(t)}\cdots a_{\textcolor{red}{j}\,i_{2}\cdots i_{m}}^{(m)}.\label{BM product}
\end{equation}
An arbitrary $m$-tuple of order $m$ hypermatrices $\left(\mathbf{A}^{(1)},\cdots,\mathbf{A}^{(m)}\right)$
for which the BM product is defined is called \emph{BM conformable}. 

We recall a variant of the BM product called the \emph{general BM
product.} The general BM product was first proposed in \cite{GER}.
It encompasses as special cases many other hypermatrix products discussed
in the literature, including the usual matrix product, the Segre outer
product, the contraction product, the higher order SVD, and the multilinear
matrix multiplication \cite{Lim2013}. In addition, the general BM
product is of particular interest to our discussion because it enables
considerable notational simplifications. The general BM product of
order $m$ hypermatrices is defined for any BM conformable $m$-tuple
$\left(\mathbf{A}^{(1)},\cdots,\mathbf{A}^{(m)}\right)$ and an additional
cubic hypermatrix $\mathbf{B}$ called the \emph{background hypermatrix}
of side length $k$ (the contracted dimension). The general BM product
denoted $\mbox{Prod}_{\mathbf{B}}\left(\mathbf{A}^{(1)},\cdots,\mathbf{A}^{(m)}\right)$
has entries given by
\[
\left[\mbox{Prod}_{\mathbf{B}}\left(\mathbf{A}^{(1)},\cdots,\mathbf{A}^{(m)}\right)\right]_{i_{1},\cdots,i_{t},\cdots,i_{m}}=
\]
\begin{equation}
\sum_{0\le\textcolor{red}{j_{1}},\textcolor{red}{\cdots},\textcolor{red}{j_{t}},\textcolor{red}{\cdots},\textcolor{red}{j_{m}}<k}a_{i_{1}\,\textcolor{red}{j_{1}}\,i_{3}\cdots i_{m}}^{(1)}\cdots\,a_{i_{1}\cdots i_{t}\,\textcolor{red}{j_{t}}\,i_{t+2}\cdots i_{m}}^{(t)}\cdots\,a_{\textcolor{red}{j_{m}}\,i_{2}\cdots i_{m}}^{(m)}b_{\textcolor{red}{j_{1}}\textcolor{red}{\cdots}\textcolor{red}{j_{t}}\textcolor{red}{\cdots}\textcolor{red}{j_{m}}}.\label{General BM product}
\end{equation}

For example, the general BM product of third order hypermatrices $\mathbf{A}^{(1)}$,
$\mathbf{A}^{(2)}$ and $\mathbf{A}^{(3)}$ with background hypermatrix
$\mathbf{B}$ denoted $\mbox{Prod}_{\mathbf{B}}\left(\mathbf{A}^{(1)},\mathbf{A}^{(2)},\mathbf{A}^{(3)}\right)$
is defined if 
\[
\mathbf{A}^{(1)}\mbox{ is}\:n_{1}\times\textcolor{red}{k}\times n_{3},\:\mathbf{A}^{(2)}\mbox{ is}\:n_{1}\times n_{2}\times\textcolor{red}{k},\:\mathbf{A}^{(3)}\mbox{ is}\:\textcolor{red}{k}\times n_{2}\times n_{3}\mbox{ and }\mbox{\ensuremath{\mathbf{B}}}\mbox{ is}\:\textcolor{red}{k}\times\textcolor{red}{k}\times\textcolor{red}{k}.
\]
The result $\mbox{Prod}_{\mathbf{B}}\left(\mathbf{A}^{(1)},\mathbf{A}^{(2)},\mathbf{A}^{(3)}\right)$
is of size $n_{1}\times n_{2}\times n_{3}$ and specified entry-wise
by 
\[
\left[\mbox{Prod}_{\mathbf{B}}\left(\mathbf{A}^{(1)},\mathbf{A}^{(2)},\mathbf{A}^{(3)}\right)\right]_{i_{1},i_{2},i_{3}}=\sum_{0\le\textcolor{red}{j_{1}},\textcolor{red}{j_{2}},\textcolor{red}{j_{3}}<k}a_{i_{1}\,\textcolor{red}{j_{1}}\,i_{3}}^{(1)}\,a_{i_{1}i_{2}\,\textcolor{red}{j_{2}}}^{(2)}\,a_{\textcolor{red}{j_{3}}\,i_{2}i_{3}}^{(3)}\,b_{\textcolor{red}{j_{1}}\textcolor{red}{j_{2}}\textcolor{red}{j_{3}}}.
\]
Note that the original BM product of order $m$ hypermatrices is recovered
from the general BM product by taking the background hypermatrix $\mathbf{B}$
to be the $m$-th order Kronecker delta hypermatrix denoted $\boldsymbol{\Delta}$,
whose entries are specified by 
\[
\left[\boldsymbol{\Delta}\right]_{i_{1},\cdots,i_{t},\cdots,i_{m}}=\begin{cases}
\begin{array}{cc}
1 & \mbox{ if }\:0\le i_{1}=\cdots=i_{t}=\cdots=i_{m}<n\\
0 & \mbox{otherwise}
\end{array}\end{cases}.
\]
In particular, Kronecker delta matrices correspond to identity matrices.

We also recall for the reader's convenience the definition of the
hypermatrix transpose operations. Let $\mathbf{A}$ be a hypermatrix
of size $n_{1}\times n_{2}\times\cdots\times n_{m}$ whose entries
are 
\[
\left[\mathbf{A}\right]_{i_{1},i_{2},\cdots,i_{m-1},i_{m}}=a_{i_{1}\,i_{2}\,\cdots\,i_{m-1}\,i_{m}}.
\]
The corresponding transpose, denoted $\mathbf{A}^{\top}$, is a hypermatrix
of size $n_{2}\times n_{3}\times\cdots\times n_{m}\times n_{1}$ whose
entries are given by 
\[
\left[\mathbf{A}^{\top}\right]_{i_{1},i_{2},\cdots,i_{m-1},i_{m}}=a_{i_{m}\,i_{1}\cdots\,i_{m-2}\,i_{m-1}}.
\]
The transpose operation  performs a cyclic permutation of the indices.
For notational convenience we adopt the convention
\[
\mathbf{A}^{\top^{2}}:=\left(\mathbf{A}^{\top}\right)^{\top},\ \mathbf{A}^{\top^{3}}:=\left(\mathbf{A}^{\top^{2}}\right)^{\top},\:\cdots,\:\mathbf{A}^{\top^{m}}:=\left(\mathbf{A}^{\top^{\left(m-1\right)}}\right)^{\top}=\mathbf{A}.
\]
By this convention 
\[
\mathbf{A}^{\top^{i}}=\mathbf{A}^{\top^{j}}\:\mbox{ if }\:i\equiv j\mod m.
\]
It follows from the definition of the transpose that
\begin{equation}
\mbox{Prod}\left(\mathbf{A}^{(1)},\mathbf{A}^{(2)},\cdots,\mathbf{A}^{(m)}\right)^{\top}=\mbox{Prod}\left(\left(\mathbf{A}^{(2)}\right)^{\top},\cdots,\left(\mathbf{A}^{(m)}\right)^{\top},\left(\mathbf{A}^{(1)}\right)^{\top}\right),\label{Transpose identity}
\end{equation}
The identity \eqref{Transpose identity} generalizes the matrix transpose
identity 
\[
\left(\mathbf{A}^{(1)}\cdot\mathbf{A}^{(2)}\right)^{\top}=\left(\mathbf{A}^{(2)}\right)^{\top}\cdot\left(\mathbf{A}^{(1)}\right)^{\top}.
\]

Finally, for notational convenience, we briefly discuss the use of
the general BM product to express multilinear forms and outer products.
Let $\mathbf{A}\in\mathbb{C}^{n_{0}\times\cdots\times n_{m-1}}$ denote
an arbitrary order $m$ hypermatrix and consider an arbitrary $m$-tuple
$\left(\mathbf{x}_{j}\in\mathbb{C}^{n_{j}\times1\times\cdots\times1}\right)_{0\le j<m}$.
The general BM product
\[
\mbox{Prod}_{\mathbf{A}}\left(\mathbf{x}_{0}^{\top^{\left(m-1\right)}},\mathbf{x}_{1}^{\top^{\left(m-2\right)}},\cdots,\mathbf{x}_{m-j-1}^{\top^{j}},\cdots,\mathbf{x}_{m-2}^{\top^{1}},\mathbf{x}_{m-1}^{\top^{0}}\right),
\]
expresses the multilinear form associated with $\mathbf{A}$. As illustration,
consider an arbitrary third order hypermatrix $\mathbf{A}\in\mathbb{C}^{m\times n\times p}$
and three vectors $\mathbf{x}\in\mathbb{C}^{m\times1\times1}$, $\mathbf{y}\in\mathbb{C}^{n\times1\times1}$
and $\mathbf{z}\in\mathbb{C}^{p\times1\times1}$ . The corresponding
multilinear form is expressed as 
\[
\mbox{Prod}_{\mathbf{A}}\left(\mathbf{x}^{\top^{2}},\mathbf{y}^{\top^{1}},\mathbf{z}^{\top^{0}}\right)=\sum_{0\le i<m}\,\sum_{0\le j<n}\,\sum_{0\le k<p}a_{ijk}\,x_{i}\,y_{j}\,z_{k}.
\]
Similarly, for an arbitrary matrix $\mathbf{A}\in\mathbb{C}^{m\times n}$
and pair of vectors $\mathbf{x}\in\mathbb{C}^{m\times1}$, $\mathbf{y}\in\mathbb{C}^{n\times1}$,
the corresponding bilinear form is expressed by 
\[
\mbox{Prod}_{\mathbf{A}}\left(\mathbf{x}^{\top^{1}},\mathbf{y}^{\top^{0}}\right)=\sum_{0\le i<m}\,\sum_{0\le j<n}a_{ij}\,x_{i}\,y_{j}=\mathbf{x}^{\top}\cdot\mathbf{A}\cdot\mathbf{y}.
\]
The general BM product also provides a convenient way to express \emph{outer
products}. For an arbitrary BM conformable $m$-tuple $\left(\mathbf{A}^{(1)},\cdots,\mathbf{A}^{(m)}\right)$,
a BM outer product corresponds to a product of the form
\begin{equation}
\mbox{Prod}\left(\mathbf{A}^{(1)}\left[:,t,:,\cdots,:\right],\,\mathbf{A}^{(2)}\left[:,:,t,\cdots,:\right],\,\cdots,\,\mathbf{A}^{(m-1)}\left[:,:,\cdots,t\right],\,\mathbf{A}^{(m)}\left[t,:,\cdots,:\right]\right).\label{Outerproduct}
\end{equation}
In \eqref{Outerproduct} we used the colon notation. Recall that in
the colon notation, $\mathbf{A}^{(1)}\left[:,t,:,\cdots,:\right]$
refers to a hypermatrix \emph{slice} of size $n_{1}\times1\times n_{3}\times\cdots\times n_{m}$
where the second index is fixed to $t$ while all other indices are
allowed to vary within their prescribed ranges. Hypermatrix outer
products are a common occurrence throughout our discussion. Fortunately,
hypermatrix outer products are conveniently expressed in terms of
general BM products. The corresponding background hypermatrices are
noted $\left\{ \boldsymbol{\Delta}^{(t)}\right\} _{0\le t<n}$ and
specified entry-wise by 
\[
\left[\boldsymbol{\Delta}^{(t)}\right]_{i_{1},\cdots,i_{m}}=\begin{cases}
\begin{array}{cc}
1 & \mbox{ if }\:0\le t=i_{1}=\cdots=i_{m}<k\\
0 & \mbox{otherwise}
\end{array}\end{cases}.
\]
The outer product in \eqref{Outerproduct} is more conveniently expressed
as $\mbox{Prod}_{\boldsymbol{\Delta}^{(t)}}\left(\mathbf{A}^{(1)},\cdots,\mathbf{A}^{(m)}\right)$.
This expression of outer products induces a natural notion of hypermatrix
BM rank. Using this notation, recall from linear algebra that a matrix
$\mathbf{B}$ is of rank $r$ (over $\mathbb{C}$) if there exists
a conformable matrix pair $\mathbf{X}^{(1)}$, $\mathbf{X}^{(2)}$
such that 
\[
\mathbf{B}=\sum_{0\le t<r}\mbox{Prod}_{\boldsymbol{\Delta}^{(t)}}\left(\mathbf{X}^{(1)},\mathbf{X}^{(2)}\right),
\]
and crucially there exists no conformable matrix pair $\mathbf{Y}^{(1)}$,
$\mathbf{Y}^{(2)}$ such that
\[
\mathbf{B}=\sum_{0\le t<r-1}\mbox{Prod}_{\boldsymbol{\Delta}^{(t)}}\left(\mathbf{Y}^{(1)},\mathbf{Y}^{(2)}\right).
\]
The definition of matrix rank above extends verbatim to hypermatrices
and is called the hypermatrix BM rank. An order $m$ hypermatrix $\mathbf{B}$
has BM rank $r$ (over $\mathbb{C}$) if there exists a BM conformable
$m$-tuple $\left(\mathbf{X}^{(1)},\cdots,\mathbf{X}^{(m)}\right)$
such that 
\[
\mathbf{B}=\sum_{0\le t<r}\mbox{Prod}_{\boldsymbol{\Delta}^{(t)}}\left(\mathbf{X}^{(1)},\cdots,\mathbf{X}^{(m)}\right),
\]
and crucially there exists no BM conformable $m$-tuple $\left(\mathbf{Y}^{(1)},\cdots,\mathbf{Y}^{(m)}\right)$
such that
\[
\mathbf{B}=\sum_{0\le t<r-1}\mbox{Prod}_{\boldsymbol{\Delta}^{(t)}}\left(\mathbf{Y}^{(1)},\cdots,\mathbf{Y}^{(m)}\right).
\]
Note that the usual notions of tensor/hypermatrix rank discussed in
the literature \cite{Lim2013} including the canonical polyadic rank
correspond to special instances of the BM rank where additional constraints
are imposed on the hypermatrices in the $m$-tuple $\left(\mathbf{X}^{(1)},\cdots,\mathbf{X}^{(m)}\right)$.

\section{General Parseval identity and Fourier transforms}

\subsection{Hypermatrix Parseval identity}

The classical matrix Parseval identity states that if $\mathbf{U}\in\mathbb{C}^{n\times n}$
is unitary then for every vectors $\mathbf{x}^{(1)}$, $\mathbf{x}^{(2)}\in\mathbb{C}^{n\times1}$ 

\[
\left(\overline{\mathbf{x}^{(1)}}\right)^{\top}\cdot\mathbf{x}^{(2)}=\left(\overline{\mathbf{U}\cdot\mathbf{x}^{(1)}}\right)^{\top}\cdot\left(\mathbf{U}\cdot\mathbf{x}^{(2)}\right).
\]

When generalizing this to hypermatrices we can't quite form the matrix-vector
products $\mathbf{U}\cdot\mathbf{x}^{(1)}$, $\mathbf{U}\cdot\mathbf{x}^{(2)}$.
Instead, notice that $\mathbf{y}^{(1)}=\mathbf{U}\cdot\mathbf{x}^{(1)}$
and $\mathbf{y}^{(2)}=\mathbf{U}\cdot\mathbf{x}^{(2)}$ satisfy

\[
\overline{y_{k}^{(1)}}\,y_{k}^{(2)}=\left[\left(\overline{\mathbf{U}\cdot\mathbf{x}^{(1)}}\right)^{\top}\right]_{k}\,\left[\mathbf{U}\cdot\mathbf{x}^{(2)}\right]_{k}=\mbox{Prod}_{\boldsymbol{\Delta}^{(k)}}\left(\left(\overline{\mathbf{U}\cdot\mathbf{x}^{(1)}}\right)^{\top},\mathbf{U}\cdot\mathbf{x}^{(2)}\right)
\]

This formulation fortunately extends to hypermatrices.

An $m$-tuple $\left(\mathbf{A}^{(1)},\cdots,\mathbf{A}^{(m)}\right)$
of order $m$ hypermatrices each cubic and of side length $n$, forms
an \emph{uncorrelated tuple} if the corresponding BM product equals
the Kronecker delta hypermatrix :
\[
\mbox{Prod}\left(\mathbf{A}^{(1)},\cdots,\mathbf{A}^{(m)}\right)=\boldsymbol{\Delta}.
\]
In some sense, uncorrelated tuples extend to hypermatrices the notion
of matrix inverse pair. Furthermore, a cubic $m$-th order hypermatrix
$\mathbf{Q}$ of side length $n$ is \emph{orthogonal} if the following
holds : 
\[
\mbox{Prod}\left(\mathbf{Q},\mathbf{Q}^{\top^{\left(m-1\right)}},\cdots,\mathbf{Q}^{\top^{k}},\cdots,\mathbf{Q}^{\top^{2}},\mathbf{Q}^{\top}\right)=\boldsymbol{\Delta}.
\]
Finally, a cubic hypermatrix $\mathbf{U}$ of even order say $2m$
and of side length $n$ is \emph{unitary} if the following holds 
\[
\mbox{Prod}\left(\mathbf{U},\overline{\mathbf{U}}^{\top^{\left(2m-1\right)}},\cdots\overline{\mathbf{U}}^{\top^{2k+1}},\mathbf{U}^{\top^{2k}},\cdots,\mathbf{U}^{\top^{2}},\overline{\mathbf{U}}^{\top}\right)=\boldsymbol{\Delta}.
\]
Both orthogonal and unitary hypermatrices yield special uncorrelated
hypermatrix tuples.\\
\\
For an arbitrary uncorrelated $m$-tuple $\left(\mathbf{A}^{(1)},\,\cdots,\mathbf{A}^{(m)}\right)$,
let $\mathbf{P}_{k}$ denote the outer product 
\[
\mathbf{P}_{k}=\mbox{Prod}_{\boldsymbol{\Delta}^{(k)}}\left(\mathbf{A}^{(1)},\cdots,\mathbf{A}^{(m)}\right),
\]
furthermore let $\left(\mathbf{x}^{(0)},\cdots,\mathbf{x}^{(m-1)}\right)$
and $\left(\mathbf{y}^{(0)},\cdots,\mathbf{y}^{(m-1)}\right)$ denote
$m$-tuples of column vectors of size $n\times1\times\cdots\times1$,
the associated\emph{ Parseval identity} is prescribed by the following
proposition\\
\textbf{}\\
\textbf{Proposition 1} : If 
\[
\forall\;0\le k<n,\quad\prod_{0\le j<m}y_{k}^{(j)}=\mbox{Prod}_{\mathbf{P}_{k}}\left(\left(\mathbf{x}^{(m-1)}\right)^{\top^{\left(m-1\right)}},\cdots,\left(\mathbf{x}^{(j)}\right)^{\top^{j}},\cdots,\left(\mathbf{x}^{(0)}\right)^{\top^{0}}\right),
\]
where $y_{k}^{(j)}$ denotes the $k$-th entry of the vector $\mathbf{y}^{(j)}$,
then we have
\[
\mbox{Prod}\left(\left(\mathbf{y}^{(m-1)}\right)^{\top^{\left(m-1\right)}},\cdots,\left(\mathbf{y}^{(j)}\right)^{\top^{j}},\cdots,\left(\mathbf{y}^{(0)}\right)^{\top^{0}}\right)
\]
\[
=
\]
\[
\mbox{Prod}\left(\left(\mathbf{x}^{(m-1)}\right)^{\top^{\left(m-1\right)}},\cdots,\left(\mathbf{x}^{(j)}\right)^{\top^{j}},\cdots,\left(\mathbf{x}^{(0)}\right)^{\top^{0}}\right).
\]
In particular, in the matrix case where $\mathbf{x}^{(0)},\mathbf{x}^{(1)}\in\mathbb{C}^{n\times1}$
and $\mathbf{A}^{(1)}$, $\mathbf{A}^{(2)}\in\mathbb{C}^{n\times n}$
are inverse pair, Parseval's identity asserts that 
\[
\forall\:\mathbf{y}^{(0)},\mathbf{y}^{(1)}\in\mathbb{C}^{n\times1}
\]
\[
\mbox{ such that }
\]
\[
\forall\:0\le k<n,\quad y_{k}^{(1)}\,y_{k}^{(0)}=\mbox{Prod}_{\mathbf{P}_{k}}\left(\left(\mathbf{x}^{(1)}\right)^{\top},\mathbf{x}^{(0)}\right)\:\mbox{ where }\mathbf{P}_{k}=\mbox{Prod}_{\boldsymbol{\Delta}^{(k)}}\left(\mathbf{A}^{\left(1\right)},\mathbf{A}^{\left(2\right)}\right)
\]
\[
\mbox{we have}
\]
\[
\mbox{Prod}\left(\left(\mathbf{y}^{(1)}\right)^{\top},\mathbf{y}^{(0)}\right)=\mbox{Prod}\left(\left(\mathbf{x}^{(1)}\right)^{\top},\mathbf{x}^{(1)}\right).
\]
\\
\\
\emph{Proof} : The proof follows from the identity
\[
\mbox{Prod}\left(\mathbf{A}^{(1)},\cdots,\mathbf{A}^{(m)}\right)=\sum_{0\le k<n}\mbox{Prod}_{\boldsymbol{\Delta}^{(k)}}\left(\mathbf{A}^{(1)},\cdots,\mathbf{A}^{(m)}\right).
\]
Consequently 
\[
\left(\sum_{0\le k<n}\,\prod_{0\le j<m}y_{k}^{(j)}\right)=\sum_{0\le k<n}\mbox{Prod}_{\mathbf{P}_{k}}\left(\left(\mathbf{x}^{(m-1)}\right)^{\top^{\left(m-1\right)}},\cdots,\left(\mathbf{x}^{(j)}\right)^{\top^{j}},\cdots,\left(\mathbf{x}^{(0)}\right)^{\top^{0}}\right)=
\]
\[
\mbox{Prod}_{\left(\sum_{0\le k<n}\mathbf{P}_{k}\right)}\left(\left(\mathbf{x}^{(m-1)}\right)^{\top^{\left(m-1\right)}},\cdots,\left(\mathbf{x}^{(j)}\right)^{\top^{j}},\cdots,\left(\mathbf{x}^{(0)}\right)^{\top^{0}}\right).
\]
This yields the desired result
\[
\mbox{Prod}\left(\left(\mathbf{y}^{(m-1)}\right)^{\top^{\left(m-1\right)}},\cdots,\left(\mathbf{y}^{(j)}\right)^{\top^{j}},\cdots,\left(\mathbf{y}^{(0)}\right)^{\top^{0}}\right)
\]
\[
=
\]
\[
\mbox{Prod}\left(\left(\mathbf{x}^{(m-1)}\right)^{\top^{\left(m-1\right)}},\cdots,\left(\mathbf{x}^{(j)}\right)^{\top^{j}},\cdots,\left(\mathbf{x}^{(0)}\right)^{\top^{0}}\right).
\]

\subsection{Hypermatrix orthogonalization and constrained uncorrelated tuples}

Applications of the proposed generalization of Parseval's identity
are predicated on the existence of non-trivial\footnote{Note that trivial orthogonal, unitary and uncorrelated hypermatrix
tuples are obtained by considering variants of Kronecker delta hypermatrices
whose nonzero entries are roots of unity.} orthogonal, unitary and uncorrelated hypermatrix tuples. We present
here an algorithmic proof of existence of non-trivial orthogonal and
uncorrelated hypermatrices of all orders and side lengths. The main
argument will be akin to proving the existence of non-trivial orthogonal
matrices by showing that the Gram\textendash Schmidt process derives
non-trivial orthogonal matrices from generic input matrices.\\
More generally, we call \emph{orthogonalization procedure}s any algorithms
which take as input some generic hypermatrices and output either orthogonal,
unitary, or uncorrelated hypermatrix tuples.

The first variant of the Gram\textendash Schmidt process which extends
to hypermatrices was proposed in \cite{phdthesis}. We will show here
that this variant of the Gram\textendash Schmidt process yields an
algorithmic proof of existence of non trivial orthogonal and non-trivial
uncorrelated hypermatrix tuples. \\
\emph{}\\
\emph{Matrix orthogonalization problem:} \\
Derive from a generic input matrix $\mathbf{A}\in\mathbb{C}^{n\times n}$
a matrix $\mathbf{X}$ of the same size subject to 
\begin{equation}
\left(\mathbf{1}_{n\times n}-\mathbf{I}_{n}\right)\circ\left(\mathbf{X}\cdot\mathbf{X}^{\top}\right)=\mathbf{0}_{n\times n},\label{Matrix_Orthogonalization}
\end{equation}
where $\circ$ denotes the entry-wise product also called the Hadamard
product, and $\mathbf{1}_{n\times n}$ denotes the $n\times n$ all
one matrix. ( Equivalently, the product Prod$\left(\mathbf{X},\,\mathbf{X}^{\top}\right)$
is a diagonal matrix. )\\
\\
\emph{Hypermatrix orthogonalization problem:}\\
Derive from a generic order $m$ input hypermatrix $\mathbf{A}\in\mathbb{C}^{n\times\cdots\times n}$
a hypermatrix $\mathbf{X}$ of the same size subject to 
\begin{equation}
\left(\mathbf{1}_{n\times\cdots\times n}-\boldsymbol{\Delta}\right)\circ\mbox{Prod}\left(\mathbf{X},\mathbf{X}^{\top^{\left(m-1\right)}},\cdots,\mathbf{X}^{\top^{2}},\mathbf{X}^{\top}\right)=\mathbf{0}_{n\times\cdots\times n},\label{Hypermatrix_Orthogonalization}
\end{equation}
where $\circ$ denotes the Hadamard product.\\
\\
It is well-known that the Gram\textendash Schmidt process yields a
solution to the matrix orthogonalization problem over $\mathbb{R}$.
We describe here a variant of the Gram\textendash Schmidt process
which extends to hypermatrices of all orders. Our proposed solution
to the matrix orthogonalization problem is obtained by solving for
the entries of $\mathbf{X}$ in the following system of $n{n \choose 2}$
equations: 
\[
\left\{ x_{u\textcolor{red}{t}}\,x_{v\textcolor{red}{t}}=a_{u\textcolor{red}{t}}\,a_{v\textcolor{red}{t}}-n^{-1}\sum_{0\le\textcolor{red}{s}<n}a_{u\textcolor{red}{s}}\,a_{v\textcolor{red}{s}}\right\} _{\begin{array}{c}
0\le\textcolor{red}{t}<n\\
0\le u<v<n
\end{array}}.
\]

(It is not hard to check that any non zero solution indeed satisfies
\eqref{Matrix_Orthogonalization}.)

For notational convenience, we rewrite the constraints above in terms
of the general BM product. The system of $n{n \choose 2}$ equations
can be more simply expressed as
\[
\forall\:0\le t<n,\;\left(\mathbf{1}_{n\times n}-\mathbf{I}_{n}\right)\circ\mbox{Prod}_{\boldsymbol{\Delta}^{(t)}}\left(\mathbf{X},\mathbf{X}^{\top}\right)=
\]
\begin{equation}
\left(\mathbf{1}_{n\times n}-\mathbf{I}_{n}\right)\circ\left[\mbox{Prod}_{\boldsymbol{\Delta}^{(t)}}\left(\mathbf{A},\mathbf{A}^{\top}\right)-n^{-1}\mbox{Prod}\left(\mathbf{A},\mathbf{A}^{\top}\right)\right],\label{Matrix orthogonalization}
\end{equation}
where $\mathbf{1}_{n\times n}$ denotes the $n\times n$ all one matrix.

Similarly, a solution to the hypermatrix orthogonalization problem
is obtained by solving for the entries of $\mathbf{X}$ in the hypermatrix
formulation of the constraints in \eqref{Matrix orthogonalization}
given by 
\[
\forall\:0\le t<n,\quad\left(\mathbf{1}_{n\times\cdots\times n}-\boldsymbol{\Delta}\right)\circ\mbox{Prod}_{\boldsymbol{\Delta}^{(t)}}\left(\mathbf{X},\mathbf{X}^{\top^{\left(m-1\right)}},\cdots,\mathbf{X}^{\top^{2}},\mathbf{X}^{\top}\right)=
\]
\begin{equation}
\left(\mathbf{1}_{n\times\cdots\times n}-\boldsymbol{\Delta}\right)\circ\left[\mbox{Prod}_{\boldsymbol{\Delta}^{(t)}}\left(\mathbf{A},\mathbf{A}^{\top^{\left(m-1\right)}},\cdots,\mathbf{A}^{\top^{2}},\mathbf{A}^{\top}\right)-n^{-1}\mbox{Prod}\left(\mathbf{A},\mathbf{A}^{\top^{\left(m-1\right)}},\cdots,\mathbf{A}^{\top^{2}},\mathbf{A}^{\top}\right)\right].\label{Hypermatrix orthogonalization}
\end{equation}

(Again, it is not hard to check that any solution to this system satisfies
\eqref{Hypermatrix_Orthogonalization}.)\\
\\
Both matrix and hypermatrix orthogonalization constraints in \eqref{Matrix orthogonalization}
and \eqref{Hypermatrix orthogonalization} turn out to be \emph{monomial
constraints}. 

General monomial constraints correspond to a system of equations which
can be expressed in terms of a coefficient matrix $\mathbf{A}\in\mathbb{C}^{m\times n}$,
a right-hand side vector $\mathbf{b}\in\mathbb{C}^{m\times1}$, and
an unknown vector $\mathbf{x}$ of size $n\times1$. These constraints
are of the form 
\begin{equation}
\left\{ \left(\prod_{0\le t<n}x_{t}^{a_{it}}\right)=b_{i}\right\} _{1\le i\le m}.\label{Monomial constraints}
\end{equation}
Such constraints are in fact linear constraints as seen by taking
the logarithm on both sides of the equal sign of each constraint.
We refer to the equivalent system obtained by taking the logarithm
as the \emph{logarithmic version} of the constraints. We solve such
systems without using logarithms to avoid any difficulty related to
branching of the logarithm. Instead, we solve such systems using a
slight variation of the Gauss\textendash Jordan elimination algorithm,
prescribed by the following elementary row operations:
\begin{itemize}
\item Row exchange : $\mbox{R}_{i}\leftrightarrow\mbox{R}_{j}$
\item Row scaling : $\left(\mbox{R}_{i}\right)^{k}\rightarrow\mbox{R}_{i}$
\item Row linear combination : $\left(\mbox{R}_{i}\right)^{k}\cdot\left(\mbox{R}_{j}\right)\rightarrow\mbox{R}_{j}$
\end{itemize}
where  $k\in\mathbb{C}$ and $\mbox{R}_{i}$ denotes the particular
constraint $\left(\prod_{0\le t<n}x_{t}^{a_{it}}\right)=b_{i}$.\\
The proposed modified row operations perform the usual row operations
prescribed by the Gauss\textendash Jordan elimination algorithm on
the logarithmic version of the constraints.\\
\\
\textbf{}\\
\textbf{Proposition 2a}: The solution $\mathbf{X}$ to the orthogonalization
constraints for a generic\footnote{A generic hypermatrix is one whose entries do not satisfy any particular
algebraic relation.} input hypermatrix $\mathbf{A}$ yields a non-trivial orthogonal hypermatrix
after normalizing of the rows of the solution matrix $\mathbf{X}$.\\
\\
\emph{Proof} : The proof follows directly from the Gaussian elimination
procedure. The row echelon form of the constraints are obtained by
performing the modified row linear combination operations described
earlier in order to put the logarithmic version of the constraints
in row echelon form. We deduce from the expression in row echelon
form of the orthogonalization constraints \eqref{Matrix orthogonalization},\eqref{Hypermatrix orthogonalization}
a criterion for establishing the existence of solutions in terms of
a single polynomial in the entries of $\mathbf{A}$ which should be
different from zero for some input. This condition will be generically
satisfied, thereby establishing the desired result.\\
\\
\\
For example, in the case of a $2\times2$ matrix 
\[
\mathbf{A}=\left(\begin{array}{rr}
a_{00} & a_{01}\\
a_{10} & a_{11}
\end{array}\right),
\]
Gauss\textendash Jordan elimination yields the solution 
\[
\mathbf{X}=\left(\begin{array}{rr}
\frac{a_{00}a_{10}-a_{01}a_{11}}{2\,x_{10}} & -\frac{a_{00}a_{10}-a_{01}a_{11}}{2\,x_{11}}\\
x_{10} & x_{11}
\end{array}\right).
\]
The rows of $\mathbf{X}$ can be normalized to form an orthogonal
matrix if no division by zero occurs and 
\[
\left[\mathbf{X}\cdot\mathbf{X}^{\top}\right]_{0,0}\ne0\,,\left[\mathbf{X}\cdot\mathbf{X}^{\top}\right]_{1,1}\ne0\Leftrightarrow\left(a_{00}a_{10}-a_{01}a_{11}\right)\left(x_{10}^{2}+x_{11}^{2}\right)x_{10}x_{11}\ne0.
\]
Similarly for a $2\times2\times2$ hypermatrix 
\[
\mathbf{A}\left[:,:,0\right]=\left(\begin{array}{rr}
a_{000} & a_{010}\\
a_{100} & a_{110}
\end{array}\right),\quad\mathbf{A}\left[:,:,1\right]=\left(\begin{array}{rr}
a_{001} & a_{011}\\
a_{101} & a_{111}
\end{array}\right),
\]
Gauss\textendash Jordan elimination yields the solution 
\[
\mathbf{X}\left[:,:,0\right]=\left(\begin{array}{rr}
\frac{\left(a_{000}a_{001}a_{100}-a_{010}a_{011}a_{110}\right)x_{101}}{a_{001}a_{100}a_{101}-a_{011}a_{110}a_{111}} & \frac{\left(a_{000}a_{001}a_{100}-a_{010}a_{011}a_{110}\right)x_{111}}{a_{001}a_{100}a_{101}-a_{011}a_{110}a_{111}}\\
\frac{a_{001}a_{100}a_{101}-a_{011}a_{110}a_{111}}{2\,x_{001}x_{101}} & -\frac{a_{001}a_{100}a_{101}-a_{011}a_{110}a_{111}}{2\,x_{011}x_{111}}
\end{array}\right),\:\mathbf{X}\left[:,:,1\right]=\left(\begin{array}{rr}
x_{001} & x_{011}\\
x_{101} & x_{111}
\end{array}\right).
\]
The rows of $\mathbf{X}$ can be normalized to form an orthogonal
matrix if no division by zero occurs and
\[
\left[\mbox{Prod}\left(\mathbf{X},\mathbf{X}^{\top^{2}},\mathbf{X}^{\top}\right)\right]_{0,0,0}\ne0,\:\left[\mbox{Prod}\left(\mathbf{X},\mathbf{X}^{\top^{2}},\mathbf{X}^{\top}\right)\right]_{1,1,1}\ne0
\]
\[
\Leftrightarrow
\]
\[
\left(a_{001}a_{100}a_{101}-a_{011}a_{110}a_{111}\right)\left(a_{000}a_{001}a_{100}-a_{010}a_{011}a_{110}\right)\left(x_{101}^{3}+x_{111}^{3}\right)\left(x_{001}x_{101}x_{011}x_{111}\right)\ne0
\]
Note that the proposed orthogonalization procedure in the matrix case
is somewhat more restrictive in comparison to the Gram\textendash Schmidt
procedure. This is seen by observing that $0\ne\det\left(\mathbf{A}\right)$
is not a sufficient condition to ensure the existence of solutions
to the orthogonalization procedure. However, the proposed orthogonalization
constraints are special instances of a more general problem called
the constrained uncorrelated tuple problem. A solution to the constrained
uncorrelated tuple problem provides a proof of existence of non-trivial
uncorrelated tuples. The constrained uncorrelated tuple problem is
specified as follows.\\
\\
\emph{Constrained inverse pair problem:} \\
Derive from a generic input matrix pair $\mathbf{A}^{(1)}$, $\mathbf{A}^{(2)}\in\mathbb{C}^{n\times n}$
matrices $\mathbf{X}^{(1)}$, $\mathbf{X}^{(2)}$ of the same size
such that
\[
\left(\mathbf{1}_{n\times n}-\mathbf{I}_{n}\right)\circ\mbox{Prod}\left(\mathbf{X}^{(1)},\mathbf{X}^{(2)}\right)=\mathbf{0}_{n\times n},
\]
\[
\mbox{which minimizes}
\]
\[
\sum_{0\le t<n}\left\Vert \left(\mathbf{1}_{n\times n}-\mathbf{I}_{n}\right)\circ\left[\mbox{Prod}_{\boldsymbol{\Delta}^{(t)}}\left(\mathbf{A}^{(1)},\mathbf{A}^{(2)}\right)-\mbox{Prod}_{\boldsymbol{\Delta}^{(t)}}\left(\mathbf{X}^{(1)},\mathbf{X}^{(2)}\right)\right]\right\Vert _{\ell_{2}}^{2}
\]
where $\circ$ denotes the Hadamard product.\\
( Equivalently, the product Prod$\left(\mathbf{X}^{(1)},\,\mathbf{X}^{(2)}\right)$
is a diagonal matrix. )\\
\\
\emph{Constrained uncorrelated tuple problem:}\\
Derive from a generic $m$-tuple of order $m$ hypermatrices $\mathbf{A}^{(1)},\cdots,\mathbf{A}^{(m)}\in$
$\mathbb{C}^{n\times\cdots\times n}$ an $m$-tulple of hypermatrices
$\left(\mathbf{X}^{(1)},\cdots,\mathbf{X}^{(m)}\right)$ of the same
sizes such that 
\[
\left(\mathbf{1}_{n\times\cdots\times n}-\boldsymbol{\Delta}\right)\circ\mbox{Prod}\left(\mathbf{X}^{(1)},\mathbf{X}^{(2)},\cdots,\mathbf{X}^{(m)}\right)=\mathbf{0}_{n\times\cdots\times n},
\]
\[
\mbox{which minimizes}
\]
\[
\sum_{0\le t<n}\left\Vert \left(\mathbf{1}_{n\times\cdots\times n}-\boldsymbol{\Delta}\right)\circ\left[\mbox{Prod}_{\boldsymbol{\Delta}^{(t)}}\left(\mathbf{A}^{(1)},\,\mathbf{A}^{(2)},\,\cdots,\mathbf{A}^{(m)}\right)-\mbox{Prod}\left(\mathbf{X}^{(1)},\,\mathbf{X}^{(2)},\,\cdots,\mathbf{X}^{(m)}\right)\right]\right\Vert _{\ell_{2}}^{2}
\]
where $\circ$ denotes the Hadamard product and the hypermatrix $\boldsymbol{\Delta}$
denotes the Kronecker delta.\\
\\
\\
\textbf{Proposition 2b} : A solution to the constrained uncorrelated
tuple problem is obtained by solving for the entries of the $m$-tuple
of hypermatrices $\left(\mathbf{X}^{(1)},\cdots,\mathbf{X}^{(m)}\right)$
in the constraints 
\[
\forall\:0\le t<n,\quad\left(\mathbf{1}_{n\times\cdots\times n}-\boldsymbol{\Delta}\right)\circ\mbox{Prod}_{\boldsymbol{\Delta}^{(t)}}\left(\mathbf{X}^{(1)},\,\mathbf{X}^{(2)},\,\cdots,\mathbf{X}^{(m)}\right)=
\]
\begin{equation}
\left(\mathbf{1}_{n\times\cdots\times n}-\boldsymbol{\Delta}\right)\circ\left[\mbox{Prod}_{\boldsymbol{\Delta}^{(t)}}\left(\mathbf{A}^{(1)},\,\mathbf{A}^{(2)},\,\cdots,\mathbf{A}^{(m)}\right)-n^{-1}\mbox{Prod}\left(\mathbf{A}^{(1)},\,\mathbf{A}^{(2)},\,\cdots,\mathbf{A}^{(m)}\right)\right].\label{Uncorrealted tuple constraints}
\end{equation}
For generic input hypermatrices $\left(\mathbf{A}^{(1)},\cdots,\mathbf{A}^{(m)}\right)$$\subset\mathbb{C}^{n\times\cdots\times n}$
the rows of $\left(\mathbf{X}^{(1)},\cdots,\mathbf{X}^{(m)}\right)$
can be normalized to obtain a non-trivial uncorrelated tuple.\\
\\
\\
\emph{Proof} : The proof again follows directly from the Gauss-Jordan
elimination procedure. The constraints in \eqref{Uncorrealted tuple constraints}
correspond to a system of $n^{m}$ \emph{monomial constraints} in
$m\cdot n^{m}$ variables. We solve such a system via Gauss-Jordan
elimination. By the argument used in the Proposition 2a we know the
hypermatrices $\left(\mathbf{X}^{(1)},\cdots,\mathbf{X}^{(m)}\right)$
can be normalized to form non-trivial uncorrelated tuples. Finally
the fact that the obtained solution minimizes the sum
\[
\sum_{0\le t<n}\left\Vert \left(\mathbf{1}_{n\times\cdots\times n}-\boldsymbol{\Delta}\right)\circ\left(\mbox{Prod}_{\boldsymbol{\Delta}^{(t)}}\left(\mathbf{A}^{(1)},\,\mathbf{A}^{(2)},\,\cdots,\mathbf{A}^{(m)}\right)-\mbox{Prod}\left(\mathbf{X}^{(1)},\,\mathbf{X}^{(2)},\,\cdots,\mathbf{X}^{(m)}\right)\right)\right\Vert _{\ell_{2}}^{2}
\]
follows from the fact that the right-hand side of equality in \eqref{Uncorrealted tuple constraints}
expresses an orthogonal projection.\\
\\
Our proposed solution to the uncorrelated tuple problem therefore
yields an algorithmic proof of existence of non trivial uncorrelated
tuples. The following corollary follows from Proposition 2a\\
\\
\textbf{Corollary 2c} : For every order $m\ge2$ and every side length
$n\ge2$ there exists an orthogonal hypermatrix having no zero entries.\\
\emph{}\\
\emph{Proof} : By Proposition 2a, The row echelon form of the constraints
\eqref{Matrix orthogonalization},\eqref{Hypermatrix orthogonalization}
yields a criterion for the existence of non trivial solution. The
criterion is expressed as a non-zero polynomial in the entries of
$\mathbf{A}$ and possibly some free variables, which in turn must
not evaluate to zero for our choice of input. A generic choice of
$\mathbf{A}$ and of free variables will indeed satisfy this requirement.

\subsection{Direct sums and Kronecker products of hypermatrices}

Recall from linear algebra that the direct sum and the Kronecker product
of square matrices\\
$\mathbf{A}\in\mathbb{C}^{n_{0}\times n_{0}}$, $\mathbf{B}\in\mathbb{C}^{n_{1}\times n_{1}}$
can both be defined in terms of bilinear forms. For notational convenience
we express here multilinear forms as general BM products. 
\[
\mbox{Prod}_{\mathbf{A}\oplus\mathbf{B}}\left(\left(\begin{array}{c}
\mathbf{x}_{1}\\
\mathbf{y}_{1}
\end{array}\right)^{\top^{1}},\left(\begin{array}{c}
\mathbf{x}_{0}\\
\mathbf{y}_{0}
\end{array}\right)^{\top^{0}}\right)\,:=\mbox{Prod}_{\mathbf{A}}\left(\mathbf{x}_{1}^{\top^{1}},\mathbf{x}_{0}^{\top^{0}}\right)+\mbox{Prod}_{\mathbf{B}}\left(\mathbf{y}_{1}^{\top^{1}},\mathbf{y}_{0}^{\top^{0}}\right)
\]
\[
\mbox{and}
\]
\[
\mbox{Prod}_{\mathbf{A}\otimes\mathbf{B}}\left(\left(\mathbf{x}_{1}\otimes\mathbf{y}_{1}\right)^{\top^{1}},\left(\mathbf{x}_{0}\otimes\mathbf{y}_{0}\right)^{\top^{0}}\right)\,:=\mbox{Prod}_{\mathbf{A}}\left(\mathbf{x}_{1}^{\top^{1}},\mathbf{x}_{0}^{\top^{0}}\right)\cdot\mbox{Prod}_{\mathbf{B}}\left(\mathbf{y}_{1}^{\top^{1}},\mathbf{y}_{0}^{\top^{0}}\right),
\]
where $\left\{ \mathbf{x}_{0},\mathbf{x}_{1}\right\} \subset\mathbb{C}^{n_{0}\times1}$
and $\left\{ \mathbf{y}_{0},\,\mathbf{y}_{1}\right\} \subset\mathbb{C}^{n_{1}\times1}$.
These definitions extend verbatim to cubic hypermatrices of all orders
as illustrated below for third order hypermatrices:
\[
\mbox{Prod}_{\mathbf{A}\oplus\mathbf{B}}\left(\left(\begin{array}{c}
\mathbf{x}_{2}\\
\mathbf{y}_{2}
\end{array}\right)^{\top^{2}},\left(\begin{array}{c}
\mathbf{x}_{1}\\
\mathbf{y}_{1}
\end{array}\right)^{\top^{1}},\left(\begin{array}{c}
\mathbf{x}_{0}\\
\mathbf{y}_{0}
\end{array}\right)^{\top^{0}}\right)\,:=\mbox{Prod}_{\mathbf{A}}\left(\mathbf{x}_{2}^{\top^{2}},\mathbf{x}_{1}^{\top^{1}},\mathbf{x}_{0}^{\top^{0}}\right)+\mbox{Prod}_{\mathbf{B}}\left(\mathbf{y}_{2}^{\top^{2}},\mathbf{y}_{1}^{\top^{1}},\mathbf{y}_{0}^{\top^{0}}\right),
\]
\[
\mbox{and}
\]
\[
\mbox{Prod}_{\mathbf{A}\otimes\mathbf{B}}\left(\left(\mathbf{x}_{2}\otimes\mathbf{y}_{2}\right)^{\top^{2}},\left(\mathbf{x}_{1}\otimes\mathbf{y}_{1}\right)^{\top^{1}},\left(\mathbf{x}_{0}\otimes\mathbf{y}_{0}\right)^{\top^{0}}\right)\,:=\mbox{Prod}_{\mathbf{A}}\left(\mathbf{x}_{2}^{\top^{2}},\mathbf{x}_{1}^{\top^{1}},\mathbf{x}_{0}^{\top^{0}}\right)\cdot\mbox{Prod}_{\mathbf{B}}\left(\mathbf{y}_{2}^{\top^{2}},\mathbf{y}_{1}^{T^{1}},\mathbf{y}_{0}^{T^{0}}\right),
\]
where $\left\{ \mathbf{x}_{0},\,\mathbf{x}_{1},\,\mathbf{x}_{2}\right\} \subset\mathbb{C}^{n_{0}\times1\times1}$
and $\left\{ \mathbf{y}_{0},\,\mathbf{y}_{1},\,\mathbf{y}_{2}\right\} \subset\mathbb{C}^{n_{1}\times1\times1}$.
\\
\\
\\
\textbf{Lemma 3} : For any two arbitrary uncorrelated $m$-tuples
of hypermatrices $\left(\mathbf{A}^{(1)},\cdots,\mathbf{A}^{(m)}\right)$
and $\left(\mathbf{B}^{(1)},\cdots,\mathbf{B}^{(m)}\right)$ the following
$m$-tuples 
\[
\left(\mathbf{A}^{(1)}\oplus\mathbf{B}^{(1)},\cdots,\mathbf{A}^{(k)}\oplus\mathbf{B}^{(k)},\cdots,\mathbf{A}^{(m)}\oplus\mathbf{B}^{(m)}\right)
\]
\[
\mbox{and}
\]
\[
\left(\mathbf{A}^{(1)}\otimes\mathbf{B}^{(1)},\cdots,\mathbf{A}^{(k)}\otimes\mathbf{B}^{(k)},\cdots,\mathbf{A}^{(m)}\otimes\mathbf{B}^{(m)}\right)
\]
also form uncorrelated hypermatrix tuples.\\
\\
\emph{Proof}: The fact that the $m$-tuple of hypermatrices $\left(\mathbf{A}^{(1)}\oplus\mathbf{B}^{(1)},\cdots,\mathbf{A}^{(k)}\oplus\mathbf{B}^{(k)},\cdots,\mathbf{A}^{(m)}\oplus\mathbf{B}^{(m)}\right)$
forms an uncorrelated tuple (assuming that the $m$-tuples $\left(\mathbf{A}^{(1)},\cdots,\mathbf{A}^{(m)}\right)$
and $\left(\mathbf{B}^{(1)},\cdots,\mathbf{B}^{(m)}\right)$ form
uncorrelated $m$-tuples) follows from the fact that the BM product
is well behaved relative to conformable block hypermatrix partitions.
Hypermatrix block partitioning schemes are hypermatrix analog of matrix
partitioning schemes into submatrices. It is convenient to think of
block partitions as hypermatrices whose entries are hypermatrices
of the same order. Let $\mathbf{U}^{(1)},\,\cdots,\mathbf{U}^{(m)}$
denote a BM conformable tuple of hypermatrices. Let $\left\{ \mathbf{U}_{i_{1}\,i_{1}\,i_{3}\,\cdots\,i_{m}}^{(t)}\right\} _{i_{1}\,i_{1}\,i_{3}\,\cdots\,i_{m}}$
denote the block partitions of the hypermatrix $\mathbf{U}^{(t)}$.
The corresponding block partition product equality is expressed by
\[
\left[\mbox{Prod}\left(\mathbf{U}^{(1)},\,\cdots,\mathbf{U}^{(m)}\right)\right]_{i_{1},\cdots,i_{m}}=\sum_{0\le\textcolor{red}{j}<k}\mbox{Prod}\left(\mathbf{U}_{i_{1}\,\textcolor{red}{j}\,i_{3}\,\cdots\,i_{m}}^{(1)}\cdots\,\mathbf{U}_{i_{1}\,\cdots\,i_{t}\,\textcolor{red}{j}\,i_{t+2}\,\cdots\,i_{m}}^{(t)}\cdots\,\mathbf{U}_{\textcolor{red}{j}\,i_{2}\,\cdots\,i_{m}}^{(m)}\right),
\]
as long as the hypermatrix blocks $\mathbf{U}_{i_{1}\,\textcolor{red}{j}\,i_{3}\,\cdots\,i_{m}}^{(1)}\cdots\,\mathbf{U}_{i_{1}\,\cdots\,i_{t}\,\textcolor{red}{j}\,i_{t+2}\,\cdots\,i_{m}}^{(t)}\cdots\,\mathbf{U}_{\textcolor{red}{j}\,i_{2}\,\cdots\,i_{m}}^{(m)}$
are always BM conformable.

Finally, the fact that the $m$-tuple of hypermatrices 
\[
\left(\mathbf{A}^{(1)}\otimes\mathbf{B}^{(1)},\cdots,\mathbf{A}^{(k)}\otimes\mathbf{B}^{(k)},\cdots,\mathbf{A}^{(m)}\otimes\mathbf{B}^{(m)}\right)
\]
also forms an uncorrelated $m$-tuple follows from the easily verifiable
BM-product identity
\[
\mbox{Prod}\left(\mathbf{A}^{(1)}\otimes\mathbf{B}^{(1)},\cdots,\mathbf{A}^{(k)}\otimes\mathbf{B}^{(k)},\cdots,\mathbf{A}^{(m)}\otimes\mathbf{B}^{(m)}\right)=
\]
\begin{equation}
\mbox{Prod}\left(\mathbf{A}^{(1)},\cdots,\mathbf{A}^{(k)},\cdots,\mathbf{A}^{(m)}\right)\otimes\mbox{Prod}\left(\mathbf{B}^{(1)},\cdots,\mathbf{B}^{(k)},\cdots,\mathbf{B}^{(m)}\right)\label{Kronecker product identity}
\end{equation}
The identity \eqref{Kronecker product identity} extends to hypermatrices
the classical matrix identity
\[
\left(\mathbf{A}^{(1)}\otimes\mathbf{B}^{(1)}\right)\cdot\left(\mathbf{A}^{(2)}\otimes\mathbf{B}^{(2)}\right)=\left(\mathbf{A}^{(1)}\cdot\mathbf{A}^{(2)}\right)\otimes\left(\mathbf{B}^{(1)}\cdot\mathbf{B}^{(2)}\right).
\]

\subsection{From matrix transformations to hypermatrix transformations}

Recall from linear algebra that we associate with some matrix $\mathbf{A}\in\mathbb{C}^{n\times n}$
a matrix transformation acting on $\mathbb{C}^{n\times1}$ defined
by the product 
\[
\forall\,\mathbf{x}\in\mathbb{C}^{n\times1},\quad\mathbf{A}\cdot\mathbf{x}
\]
In order to extend to hypermatrices the notion of transformation acting
on a vector space, we reformulate the matrix transformations above
as follows:
\[
\mathcal{T}_{\mathbf{A}^{\top},\mathbf{A}}\,:\,\mathbb{C}^{n\times1}\rightarrow\mathbb{C}^{n\times1}
\]
\[
\mathbf{y}=\mathcal{T}_{\mathbf{A}^{\top},\mathbf{A}}\left(\mathbf{x}\right)\Leftrightarrow\forall\;0\le k<n,\;y_{k}=\sqrt{\mbox{Prod}_{\mathbf{P}_{k}}\left(\mathbf{x}^{\top},\mathbf{x}\right)},
\]
where $\mathbf{P}_{k}=\mbox{Prod}_{\boldsymbol{\Delta}^{(k)}}\left(\mathbf{A}^{\top},\mathbf{A}\right)$.
Consequently, up to sign,
\[
\mathbf{y}=\mathbf{A}\cdot\mathbf{x}.
\]

(That is, this equation holds if we identify two complex numbers differing
only by sign.)

Note that linear transformations such as $\mathcal{T}_{\mathbf{A}^{\top},\mathbf{A}}$
are special cases of an equivalence classes of non-linear transformations
associated with an arbitrary pair of $n\times n$ matrices $\mathbf{A}^{(1)}$,
$\mathbf{A}^{(2)}$ defined by 
\[
\mathcal{T}_{\mathbf{A}^{(1)},\mathbf{A}^{(2)}}\,:\,\mathbb{C}^{n\times1}\rightarrow\mathbb{C}^{n\times1}
\]
\[
\mathbf{y}=\mathcal{T}_{\mathbf{A}^{(1)},\mathbf{A}^{(2)}}\left(\mathbf{x}\right)\Leftrightarrow\forall\;0\le k<n,\quad y_{k}=\sqrt{\mbox{Prod}_{\mathbf{P}_{k}}\left(\mathbf{x}^{\top},\mathbf{x}\right)},
\]
where $\mathbf{P}_{k}=\mbox{Prod}_{\boldsymbol{\Delta}^{(k)}}\left(\mathbf{A}^{(1)},\mathbf{A}^{(2)}\right)$.
Such equivalence classes of transformations naturally extend to hypermatrices
and are motivated by the general Parseval identity. We define for
an arbitrary $m$-tuple of order $m$ hypermatrices $\left(\mathbf{A}^{(1)},\cdots,\mathbf{A}^{(m)}\right)$
the equivalence class of transforms $\mathcal{T}_{\mathbf{A}^{(1)},\cdots,\mathbf{A}^{(m)}}$
whose action on the vector space $\mathbb{C}^{n\times1\times\cdots\times1}$
is defined by 
\[
\mathcal{T}_{\mathbf{A}^{(1)},\cdots,\mathbf{A}^{(m)}}\,:\,\mathbb{C}^{n\times1\times\cdots\times1}\rightarrow\mathbb{C}^{n\times1\times\cdots\times1}
\]
\[
\mbox{such that }
\]
\[
\mathbf{y}=\mathcal{T}_{\mathbf{A}^{(1)},\cdots,\mathbf{A}^{(m)}}\left(\mathbf{x}\right)
\]
\[
\Leftrightarrow
\]
\[
\forall\;0\le k<n,\quad y_{k}=\sqrt[m]{\mbox{Prod}_{\mathbf{P}_{k}}\left(\mathbf{x}^{\top^{\left(m-1\right)}},\mathbf{x}^{\top^{\left(m-2\right)}},\cdots,\mathbf{x}^{\top^{j}},\cdots,\mathbf{x}^{\top^{1}},\mathbf{x}^{\top^{0}}\right)}
\]
where $\mathbf{P}_{k}=\mbox{Prod}_{\boldsymbol{\Delta}^{(k)}}\left(\mathbf{A}^{(1)},\cdots,\mathbf{A}^{(m)}\right)$.
The equivalence class of transforms associated with $m$-th order
hypermatrices is defined modulo multiplication of the each entry of
the image vector $\mathbf{y}$ with an arbitrary $m$-th root of unity.

\subsection{Hypermatrix Fourier transforms}

Hypermatrix transforms also motivate a natural generalization of Fourier
transforms. To emphasize the analogy between the hypermatrix Fourier
transform and the matrix Fourier transform we briefly recall here
a matrix variant of the Fourier transform. Given an inverse pair of
$n\times n$ matrices (i.e. an uncorrelated pair of second order hypermatrices
$\left(\mathbf{A}^{(1)},\mathbf{A}^{(2)}\right)$) their induced Fourier
transform, denoted $\mathcal{T}_{\mathbf{A}^{(1)},\mathbf{A}^{(2)}}$,
is defined as the map acting on the vector space $\mathbb{C}^{n\times1}$
defined by
\[
\mathcal{T}_{\mathbf{A}^{(1)},\mathbf{A}^{(2)}}\,:\,\mathbb{C}^{n\times1}\rightarrow\mathbb{C}^{n\times1}
\]
\[
\mbox{such that }
\]
\[
\mathbf{y}=\mathcal{T}_{\mathbf{A}^{(1)},\mathbf{A}^{(2)}}\left(\mathbf{x}\right)\Leftrightarrow\forall\;0\le k<n,\quad y_{k}=\sqrt{\mbox{Prod}_{\mathbf{P}_{k}}\left(\mathbf{x}^{\top},\mathbf{x}\right)},
\]
where $\mathbf{P}_{k}=\mbox{Prod}_{\boldsymbol{\Delta}^{(k)}}\left(\mathbf{A}^{(1)},\mathbf{A}^{(2)}\right)$.
Although different choices of branches for the square root induce
different transforms we consider all such transforms to belong to
the same equivalence class of transforms for which 
\[
\forall\;0\le k<n,\;\left(y_{k}\right)^{2}=\mbox{Prod}_{\mathbf{P}_{k}}\left(\mathbf{x}^{\top},\mathbf{x}\right).
\]
In linear algebra terms, we say that such maps are equivalent up to
multiplication of the image vector $\mathbf{y}$ by a diagonal matrix
whose diagonal entries are either $-1$ or $1$. Furthermore, by Parseval's
identity we know that the transform $\mathcal{T}_{\mathbf{A}^{(1)},\mathbf{A}^{(2)}}$
preserves the sum of squares of entries of the pre-image $\mathbf{x}$
: 
\[
\mathbf{y}=\mathcal{T}_{\mathbf{A}^{(1)},\mathbf{A}^{(2)}}\left(\mathbf{x}\right)\Leftrightarrow\mbox{Prod}\left(\mathbf{y}^{\top},\mathbf{y}\right)=\mbox{Prod}\left(\mathbf{x}^{\top},\mathbf{x}\right)
\]
Similarly, we associate with some arbitrary uncorrelated $m$-tuples
of hypermatrices $\left(\mathbf{A}^{(1)},\cdots,\mathbf{A}^{(m)}\right)$,
each of order $m$ and having side length $n$, a \emph{hypermatrix
Fourier transform} denoted $\mathcal{T}_{\mathbf{A}^{(1)},\cdots,\mathbf{A}^{(m)}}$
whose action on the vector space $\mathbb{C}^{n\times1\times\cdots\times1}$
is defined by 
\[
\mathcal{T}_{\mathbf{A}^{(1)},\cdots,\mathbf{A}^{(m)}}\,:\,\mathbb{C}^{n\times1\times\cdots\times1}\rightarrow\mathbb{C}^{n\times1\times\cdots\times1}
\]
\[
\mbox{such that }
\]
\[
\mathbf{y}=\mathcal{T}_{\mathbf{A}^{(1)},\cdots,\mathbf{A}^{(m)}}\left(\mathbf{x}\right)
\]
\[
\Leftrightarrow
\]
\[
\forall\;0\le k<n,\quad y_{k}=\sqrt[m]{\mbox{Prod}_{\mathbf{P}_{k}}\left(\mathbf{x}^{\top^{\left(m-1\right)}},\mathbf{x}^{\top^{\left(m-2\right)}},\cdots,\mathbf{x}^{\top^{j}},\cdots,\mathbf{x}^{\top^{1}},\mathbf{x}^{\top^{0}}\right)},
\]
where $\mathbf{P}_{k}=\mbox{Prod}_{\boldsymbol{\Delta}^{(k)}}\left(\mathbf{A}^{(1)},\cdots,\mathbf{A}^{(m)}\right)$.
Although different choices of branches for the $m$-th root induce
different transforms, we consider such transforms to belong to the
equivalence class of transforms for which 
\[
\forall\;0\le k<n,\quad\left(y_{k}\right)^{m}=\mbox{Prod}_{\mathbf{P}_{k}}\left(\mathbf{x}^{\top^{\left(m-1\right)}},\mathbf{x}^{\top^{\left(m-2\right)}},\cdots,\mathbf{x}^{\top^{j}},\cdots,\mathbf{x}^{\top^{1}},\mathbf{x}^{\top^{0}}\right).
\]
These transforms are equivalent up to multiplication of each entry
of the image $\mathbf{y}$ by an arbitrary $m$-th root of unity.
By Proposition 1 it follows that the proposed transform preserves
the sum of $m$-th powers of entries of $\mathbf{x}$:
\[
\mathbf{y}=\mathcal{T}_{\mathbf{A}^{(1)},\cdots,\mathbf{A}^{(m)}}\left(\mathbf{x}\right)
\]
\[
\Leftrightarrow
\]
\[
\mbox{Prod}\left(\mathbf{x}^{\top^{\left(m-1\right)}},\mathbf{x}^{\top^{\left(m-2\right)}},\cdots,\mathbf{x}^{\top^{j}},\cdots,\mathbf{x}^{\top^{1}},\mathbf{x}^{\top^{0}}\right)=\mbox{Prod}\left(\mathbf{y}^{\top^{\left(m-1\right)}},\mathbf{y}^{\top^{\left(m-2\right)}},\cdots,\mathbf{y}^{\top^{j}},\cdots,\mathbf{y}^{\top^{1}},\mathbf{y}^{\top^{0}}\right).
\]

\subsection{Third order DFT hypermatrices}

We recall from matrix algebra that matrix inverse pairs associated
with the Discrete Fourier Transform (DFT) acting on the vector space
$\mathbb{C}^{n\times1}$ corresponds to $\mathcal{T}_{\mathbf{F},\overline{\mathbf{F}}^{\top}}$
where the entries of the $n\times n$ matrix $\mathbf{F}$ are given
by 
\[
\left[\mbox{\ensuremath{\mathbf{F}}}\right]_{u,v}=\frac{1}{\sqrt{n}}\exp\left\{ i\,\frac{2\pi}{n}\,u\,v\right\} .
\]
The definition crucially relies on the following geometric sum identity
valid for every non zero integer $n$ 
\[
\left(\frac{1}{n}\sum_{0\le t<n}\exp\left\{ i\,\frac{2\pi}{n}\,u\,t-i\,\frac{2\pi}{n}\,t\,v\right\} \right)=\begin{cases}
\begin{array}{cc}
1 & \mbox{ if }\:0\le u=v<n\\
0 & \mbox{otherwise}
\end{array}.\end{cases}
\]
\[
\Leftrightarrow
\]
\begin{equation}
\sum_{0\le t<n}\left(\frac{\exp\left\{ i\,\frac{2\pi}{n}\,u\,t\right\} }{\sqrt{n}}\right)\,\left(\frac{\exp\left\{ -i\,\frac{2\pi}{n}\,t\,v\right\} }{\sqrt{n}}\right)=\begin{cases}
\begin{array}{cc}
1 & \mbox{ if }\:0\le u=v<n\\
0 & \mbox{otherwise}
\end{array}.\end{cases}\label{Matrix DFT identity}
\end{equation}
The equality above expresses the fact that the $n\times n$ matrices
$\mathbf{F}$ and $\overline{\mathbf{F}}^{\top}$ are in fact inverse
pairs (i.e. uncorrelated pair of second order hypermatrices). We therefore
understand the DFT to be associated with a special Fourier transform.
In this Fourier transform the entries of the inverse matrix pairs
are roots of unity scaled by the normalizing factor $\nicefrac{1}{\sqrt{n}}$.
By Lemma 3 if $\mathcal{T}_{\mathbf{F},\overline{\mathbf{F}}^{\top}}$
is a DFT then for every integer $k>1$ the Fourier transform $\mathcal{T}_{\mathbf{F}^{\otimes^{k}},\overline{\mathbf{F}^{\otimes^{k}}}^{\top}}$
is a also a DFT. Recall that $\otimes^{k}$ means $k$ repeated Kronecker
products. 

There is a third order hypermatrix identity similar to the identity
in \eqref{Matrix DFT identity}, which is valid for values of the
positive integer $n$ characterized in Proposition 4. The third order
DFT hypermatrix identity crucially relies on the following geometric
sum identity
\[
\left(\frac{1}{n}\sum_{0\le t<n}\exp\left\{ i\frac{2\pi}{n}\left(u\,\sqrt{t}-\sqrt{t}\,w\right)^{2}+i\frac{2\pi}{n}\left(u\,\sqrt{t}-v\,\sqrt{t}\right)^{2}+i\frac{2\pi}{n}\left(\sqrt{t}\,v-\sqrt{t}\,w\right)^{2}\right\} \right)=
\]
\[
\begin{cases}
\begin{array}{cc}
1 & \mbox{ if }\:0\le u=v=w<n\\
0 & \mbox{otherwise}
\end{array},\end{cases}
\]
for values of $n$ characterized in Proposition 4. The identity above
can be rewritten as
\[
\left(\sum_{0\le t<n}\frac{\exp\left\{ i\frac{2\pi}{n}\left(u\,\sqrt{t}-\sqrt{t}\,w\right)^{2}\right\} }{\sqrt[3]{n}}\,\frac{\exp\left\{ i\frac{2\pi}{n}\left(u\,\sqrt{t}-v\,\sqrt{t}\right)^{2}\right\} }{\sqrt[3]{n}}\,\frac{\exp\left\{ i\frac{2\pi}{n}\left(\sqrt{t}\,v-\sqrt{t}\,w\right)^{2}\right\} }{\sqrt[3]{n}}\right)=
\]
\[
\begin{cases}
\begin{array}{cc}
1 & \mbox{ if }\:0\le u=v=w<n\\
0 & \mbox{otherwise}
\end{array}.\end{cases}
\]
The identity above expresses a BM product of the uncorrelated triple
$\left(\mathbf{F},\mathbf{G},\mathbf{H}\right)$. Note that the entries
of $\mathbf{F}$, $\mathbf{G}$ and $\mathbf{H}$ are $n$-th roots
of unity scaled by the same normalizing factor $\nicefrac{1}{\sqrt[3]{n}}$.
The entries of $\mathbf{F}$, $\mathbf{G}$ and $\mathbf{H}$ are
thus given by 
\begin{equation}
\left[\mathbf{F}\right]_{u,t,w}=\frac{\exp\left\{ i\frac{2\pi}{n}t\,\left(u-w\right)^{2}\right\} }{\sqrt[3]{n}},\;\left[\mathbf{G}\right]_{u,v,t}=\frac{\exp\left\{ i\frac{2\pi}{n}t\,\left(u-v\right)^{2}\right\} }{\sqrt[3]{n}},\:\left[\mathbf{H}\right]_{t,v,w}=\frac{\exp\left\{ i\frac{2\pi}{n}t\left(v-w\right)^{2}\right\} }{\sqrt[3]{n}}.\label{Third order DFT}
\end{equation}
As a result the transform $\mathcal{T}_{\mathbf{F},\mathbf{G},\mathbf{H}}$
is a hypermatrix DFT acting on the vector space $\mathbb{C}^{n\times1\times1}$.
The smallest possible choice for $n$ is $n=5$. By Lemma 3, if $\mathcal{T}_{\mathbf{F},\mathbf{G},\mathbf{H}}$
is a DFT over $\mathbb{C}^{n\times1\times1}$ then for every positive
integer $k>1$, $\mathcal{T}_{\mathbf{F}^{\otimes^{k}},\mathbf{G}^{\otimes^{k}},\mathbf{H}^{\otimes^{k}}}$
is also a DFT over the vector space $\mathbb{C}^{n^{k}\times1\times1}$.
\\
The following proposition determines the necessary and sufficient
condition on the positive integer $n$ which ensures that the hypermatrices
in \eqref{Third order DFT} are uncorrelated.\\

\textbf{Proposition 4} : The $n\times n\times n$ hypermatrices $\mathbf{F}$,
$\mathbf{G}$ and $\mathbf{H}$ whose entries are specified by 
\[
\left[\mathbf{F}\right]_{u,t,w}=\frac{\exp\left\{ i\frac{2\pi}{n}t\,\left(u-w\right)^{2}\right\} }{\sqrt[3]{n}},\;\left[\mathbf{G}\right]_{u,v,t}=\frac{\exp\left\{ i\frac{2\pi}{n}t\,\left(u-v\right)^{2}\right\} }{\sqrt[3]{n}},\:\left[\mathbf{H}\right]_{t,v,w}=\frac{\exp\left\{ i\frac{2\pi}{n}t\,\left(v-w\right)^{2}\right\} }{\sqrt[3]{n}}.
\]
form an uncorrelated triple if and only if the equation 
\[
x^{2}+3\,y^{2}\equiv0\mod n,
\]
admits no solution other then the trivial solution $x\equiv0\mod n$
and $y\equiv0\mod n$.\\

\emph{Proof}: The construction requires the following implication
\[
\forall\;0\le u,v,w<n,\,u\left(e^{i\frac{2\pi}{3}}\right)^{2}+v\left(e^{i\frac{2\pi}{3}}\right)^{1}+w\left(e^{i\frac{2\pi}{3}}\right)^{0}\ne0\Rightarrow\left(u-v\right)^{2}+\left(v-w\right)^{2}+\left(u-w\right)^{2}\ne0\mbox{ mod }n.
\]
Let $x=u-v$ and $y=v-w$, the implication becomes 
\[
\forall\;x,y\in\mathbb{N},\quad x^{2}+y^{2}+\left(x+y\right)^{2}\ne0\mod n.
\]
\[
\Rightarrow\forall\;x,y\in\mathbb{N},\quad2\left(x^{2}+xy+y^{2}\right)\ne0\mod n.
\]
If $n$ is even then the choice $x=\frac{n}{2}$ and $y=0$ always
constitutes a counterexample. However if $n$ is odd the constraints
may be stated as follows: 

For all integers $x$, $y$ not both zero modulo $n$ we require that
\[
x^{2}+xy+y^{2}\ne0\mod n.
\]
\[
\Rightarrow\left(x+\frac{y}{2}\right)^{2}+3\,\left(\frac{y}{2}\right)^{2}\ne0\mod n.
\]
from which the sought after result follows.

In particular, when $n$ is prime we need $-3$ to be a quadratic
non-residue modulo $n$. An easy calculation shows that the primes
of the forms $12m+5$ and $12m+11$ satisfy these conditions, and
in particular there are infinitely many such $n$. We leave the case
of composite $n$ to the reader.

\subsection{Hadamard hypermatrices}

We discuss here Hadamard hypermatrices which are used to construct
special DFT hypermatrices which have real entries. In fact we extend
to hypermatrices Sylvester's classical Hadamard matrix construction.
Recall from linear algebra that a matrix $\mathbf{H}\in\left\{ -1,1\right\} ^{n\times n}$
is a Hadamard matrix if 
\begin{equation}
\left[\mathbf{H}\cdot\mathbf{H}^{\top}\right]_{i,j}=\begin{cases}
\begin{array}{cc}
n & \mbox{ if }\:0\le i=j<n\\
0 & \mbox{otherwise}
\end{array}.\end{cases}\label{Hadamard constraint}
\end{equation}
Hadamard matrices are of considerable importance in topics relating
to combinatorial design and the analysis of boolean functions. They
are also used to define the famous Hadamard\textendash Rademacher\textendash Walsh
transform which plays an important role in Quantum computing and signal
processing. Hadamard matrices are also common occurrences in practical
implementations of the Fast Fourier Transform. Furthermore, Hadamard
matrices are well-known to be optimal matrices relative to the Hadamard
determinant inequality 
\[
\left|\det\boldsymbol{\Theta}\right|\le\left(\sqrt{n}\right)^{n},
\]
valid over the set of all $n\times n$ matrices $\boldsymbol{\Theta}$
whose entries of are bounded in absolute value by $1$. Equality is
achieved in Hadamard's determinant inequality for Hadamard matrices.
In 1867, James Joseph Sylvester proposed the classical construction
of an infinite family of Hadamard matrices of size $2^{n}\times2^{n}$
for any integer $n\ge1$. Sylvester's construction starts with the
$2\times2$ matrix
\[
\left(\begin{array}{cc}
1 & 1\\
1 & -1
\end{array}\right)
\]
and considers the sequence of matrices 
\[
\left\{ \left(\begin{array}{cc}
1 & 1\\
1 & -1
\end{array}\right)^{\otimes^{n}}\in\left\{ -1,1\right\} ^{2^{n}\times2^{n}}\right\} _{1<n<\infty}.
\]
By Lemma 3 we know that every matrix in the sequence above will satisfy
the Hadamard criterion \eqref{Hadamard constraint}. Having defined
in section 3.1 orthogonal hypermatrices, it is relatively straightforward
to extend the Hadamard criterion \eqref{Hadamard constraint} to hypermatrices
of arbitrary orders, which can be used to extend to hypermatrices
the Hadamard\textendash Rademacher\textendash Walsh transform. Formally,
an order $m$ hypermatrix $\mathbf{H}\in\left\{ -1,1\right\} ^{n\times\cdots\times n}$
is Hadamard if 
\begin{equation}
\left[\mbox{Prod}\left(\mathbf{H},\mathbf{H}^{\top^{\left(m-1\right)}},\cdots,\mathbf{H}^{\top^{k}},\cdots,\mathbf{H}^{\top^{2}},\mathbf{H}^{\top}\right)\right]_{i_{1},\cdots,i_{m}}=\begin{cases}
\begin{array}{cc}
n & \mbox{ if }\:0\le i_{1}=\cdots=i_{m}<n\\
0 & \mbox{otherwise}
\end{array}.\end{cases}\label{Hadamard criterion}
\end{equation}
The following theorem extends the scope of both Sylvester's constructions
and the famous Hadamard matrix conjecture.\\
\\
\textbf{}\\
\textbf{Theorem 5} : For every positive integer $n\ge1$ and every
positive integer $m$ which is either odd or equal to $2$, there
exists an order $m$ Hadamard hypermatrix of side length $2^{n}$.
In contrast, if $m$ is an even integer larger than $2$, then there
is no order $m$ Hadamard hypermatrix of side length $2$.\\
\\
\\
\emph{Proof} : By Lemma 3, it suffices to provide an explicit construction
for odd order Hadamard hypermatrices of side length $2$. For side
length $2$ hypermatrices of order $m>2$, the Hadamard criterion
\eqref{Hadamard criterion} is expressed as follows
\[
\forall\left(i_{1},\cdots,i_{m}\right)\notin\left\{ \left(0,0,\cdots,0\right),\,\left(1,1,\cdots,1\right)\right\} 
\]
\[
\left(h_{i_{1}\,0\,i_{3}\cdots i_{m}}\,h_{i_{2}\,0\,i_{4}\cdots i_{m}i_{1}}\cdots h_{i_{m}\,0\,i_{2}\cdots i_{m-1}}\right)+\left(h_{i_{1}\,1\,i_{3}\cdots i_{m}}\,h_{i_{2}\,1\,i_{4}\cdots i_{m}i_{1}}\cdots h_{i_{m}\,1\,i_{2}\cdots i_{m-1}}\right)=0
\]
\[
\mbox{and}
\]
\[
\forall\:0\le i<2,\quad\left(h_{i\,0\,i\cdots i\,i}\right)^{m}+\left(h_{i\,1\,i\cdots i\,i}\right)^{m}=2.
\]

The first set of constraints are equivalently expressed as

\[
\forall\left(i_{1},\cdots,i_{m}\right)\notin\left\{ \left(0,0,\cdots,0\right),\,\left(1,1,\cdots,1\right)\right\} 
\]
\[
\nicefrac{\left(h_{i_{1}\,0\,i_{3}\cdots i_{m}}\,h_{i_{2}\,0\,i_{4}\cdots i_{m}i_{1}}\cdots h_{i_{m}\,0\,i_{2}\cdots i_{m-1}}\right)}{\left(h_{i_{1}\,1\,i_{3}\cdots i_{m}}\,h_{i_{2}\,1\,i_{4}\cdots i_{m}i_{1}}\cdots h_{i_{m}\,1\,i_{2}\cdots i_{m-1}}\right)}=-1
\]

For $\pm1$ solutions, the second set of constraints just states that

\[
\forall\:0\le i<2,\quad\left(h_{i\,0\,i\cdots i\,i}\right)^{m}=\left(h_{i\,1\,i\cdots i\,i}\right)^{m}=1.
\]

For all $j_{1},j_{2},\ldots,j_{m-1},$ define $H_{j_{1}\,j_{2}\cdots j_{m-1}}=\nicefrac{h_{j_{1}\,0\,j_{2}\cdots j_{m-1}}}{h_{j_{1}\,1\,j_{2}\cdots j_{m-1}}}$.
The first set of constraints simplifies to

\[
H_{i_{1}i_{2}\cdots i_{m-1}}H_{i_{2}i_{3}\cdots i_{m}}\cdots H_{i_{m}i_{1}\cdots i_{m-2}}=-1\qquad\forall\left(i_{1},\cdots,i_{m}\right)\notin\left\{ \left(0,0,\cdots,0\right),\,\left(1,1,\cdots,1\right)\right\} 
\]

The second set of constraints states that

\[
\forall\:0\le i<2,\quad\left(H_{i\,i\cdots i\,i}\right)^{m}=1.
\]

Clearly the original constraints (in the original variables $h$)
have a $\pm1$ solution if and only if the new constraints (in the
new variables $H$) have a $\pm1$ solution.

We now show that if $m>2$ is even then there are no solutions. Let
$m=2k$, and consider the constraint corresponding to $i_{1}=1,i_{2}=\cdots=i_{k}=0$,$i_{k+1}=1,i_{k+2}=\cdots=i_{m}=0$.
This constraint states that

\[
H_{10\cdots0}^{2}H_{0\cdots0}^{2}H_{0\cdots01}^{2}\cdots H_{010\cdots0}^{2}=-1,
\]

which clearly has no $\pm1$ solution.

From now on, assume that $m>1$ is odd. We immediately get that

\[
H_{0\cdots0}=H_{1\cdots1}=1.
\]

Let us call a binary word of length $m$ a \emph{necklace} if it is
lexicographically smaller than all its rotations. Since rotations
of a word $i_{1}\cdots i_{m}$ correspond to the same constraint,
it is enough to consider constraints corresponding to necklaces. For
each word $i_{1}\cdots i_{m}$, a \emph{window} consists of $m-1$
contiguous characters (where contiguity is cyclic). Thus there are
$m$ windows, some of which could be the same.

If a necklace is periodic with minimal period $p$, then each window
will appear (at least) $\nicefrac{m}{p}$ times. The following lemma
shows that periodicity is the only reason that a window repeats.\\

\textbf{Lemma 5a} : Suppose that a word $w_{0}\ldots w_{m-1}$ satisfies
$w_{i}=w_{i-p}$ for $i=1,\ldots,m-1$ (but not necessarily for $i=0$).
Then $w$ has a period $\pi$ (possibly $m$) such that $p$ is a
multiple of $\pi$.\\
\\
\emph{Proof of Lemma 5a} : The proof is by induction on $m$. We can
assume $0\leq p<m$. If $m=1$ then there is nothing to prove. If
$p$ divides $m$ then the constraints imply that $p$ is a period
of $w$, so again there is nothing to prove. Suppose therefore that
$q=m\mod p>0$. The constraints imply that $w$ has the form

\[
w_{0}w_{1}\ldots w_{p-1}w_{0}w_{1}\ldots w_{p-1}\cdots w_{0}w_{1}\ldots w_{q-1},
\]

and furthermore $w_{1}=w_{q+1},\ldots,w_{p-1}=w_{p-1+q\mod p}$. That
is, the word $w_{0}\ldots w_{p-1}$ satisfies the premise of the lemma
with the shift $q$. By induction, $w_{0}\ldots w_{p-1}$ has period
$\pi$ (which thus divides $p$) and $q$ is a multiple of $\pi$.
It follows that $\pi$ divides $m$ and so is a period of $w_{0},\ldots,w_{m-1}$.
This completes the proof of the lemma.\\
\\
The lemma 5a implies that indeed if a necklace has minimal period
$p$ (possibly $p=m$) then each window appears $\nicefrac{m}{p}$
times, and so an odd number of times. We can thus restate the constraints
as follows, for $\pm1$ solutions:

\emph{For each non-constant necklace $i_{1}\ldots i_{m}$, the product
of $H$-values corresponding to distinct non-constant windows of $i_{1}\ldots i_{m}$
equals $-1$.}

As an example, for $m=5$ the non-constant necklaces are $00001,00011,00101,00111,01011,01111$,
and the corresponding constraints are

\[
H_{0001}H_{0010}H_{0100}H_{1000}=-1
\]

\[
H_{0001}H_{0011}H_{0110}H_{1100}H_{1000}=-1
\]

\[
H_{0010}H_{0101}H_{1010}H_{0100}H_{1001}=-1
\]

\[
H_{0011}H_{0111}H_{1110}H_{1100}H_{1001}=-1
\]

\[
H_{0101}H_{1011}H_{0110}H_{1101}H_{1010}=-1
\]

\[
H_{0111}H_{1110}H_{1101}H_{1011}=-1
\]

Consider now the graph whose vertex set consists of all non-constant
necklaces, and edges connect two necklaces $x,y$ if some rotations
of $x,y$ have Hamming distance 1. For example, $00101$ and $01011$
are connected since $01010$ and $01011$ differ in only one position.
It is not hard to check that each non-constant window appears in exactly
two constraints (corresponding to its two completions), and these
constraints correspond to an edge. Continuing our example, the window
$0101$ appears in the constraints corresponding to the necklaces
$00101,\,01011$, and only there. (An edge can correspond to several
windows: for example $\left(00001,\,00011\right)$ corresponds to
both $0001$ and $1000$). We will show that this graph contains a
sub-graph in which all degrees are odd. If we set to $-1$ all variables
corresponding to the chosen edges (one window per edge) and set to
$1$ all the other variables, then we obtain a solution to the set
of constraints.

For example, the edges $\left\{ \left(00001,\,00011\right),\,\left(00101,\,00111\right),\,\left(01011,\,01111\right)\right\} $
constitute a matching in the graph, and so setting $H_{0001}=H_{1001}=H_{1011}=-1$
and setting all other variables to $1$ yields a solution.

A well-known result states that a connected graph contains a sub-graph
in which all degrees are odd if and only if it has an even number
of vertices\footnote{Here is a proof of the hard direction, taken from Jukna's \emph{Extremal
combinatorics}: partition the graph into a list of pairs, and choose
a path connecting each pair. Now take the XOR of all these paths.}. To complete the proof, it thus suffices to show that the number
of necklaces (and so non-constant necklaces) is even. The classical
formula for the number of necklaces (obtainable using the orbit-stabilizer
theorem) states that the number of binary necklaces of length $m$
is

\[
\frac{1}{m}\sum_{k|m}\varphi(k)\,2^{\frac{m}{k}}.
\]

Here $\varphi$ is Euler's function. Since $m$ is odd, it suffices
to show that all summands are even. This is clear for all summands
with $k<m$. When $k=m$, we use the easy fact that $\varphi(m)$
is even for all $m>2$, which follows from the explicit formula for
$\varphi(m)$ in terms of the factorization of $m$. This completes
the proof of Theorem 5.\\
\\
We close this section with an explicit example of a $2\times2\times2$
Hadamard hypermatrix: 
\[
\mathbf{H}\left[:,:,0\right]=\left(\begin{array}{cc}
1 & 1\\
-1 & 1
\end{array}\right),\quad\mathbf{H}\left[:,:,1\right]=\left(\begin{array}{cc}
1 & 1\\
1 & 1
\end{array}\right).
\]

\section{Spectral decomposition of Kronecker products and direct sums of side
length $2$ hypermatrices}

We describe here elementary methods for deriving generators for matrix
and hypermatrix spectral elimination ideals, which will be defined
here.

\subsection{The matrix case}

We start by describing the derivation of generators for the matrix
\emph{spectral elimination ideal} which we now define. Let $\mathbf{A}\in\mathbb{C}^{2\times2}$
having distinct eigenvalues $\lambda_{0}$, $\lambda_{1}$. For the
purposes of our derivation the eigenvalues will be expressed as
\[
\lambda_{0}=\mu_{0}\cdot\nu_{0},\quad\lambda_{1}=\mu_{1}\cdot\nu_{1}.
\]
Recall that the spectral decomposition equation is given by 
\begin{equation}
\left(\begin{array}{cc}
a_{00} & a_{01}\\
a_{10} & a_{11}
\end{array}\right)=\left[\left(\begin{array}{cc}
u_{00} & u_{01}\\
u_{10} & u_{11}
\end{array}\right)\cdot\left(\begin{array}{cc}
\mu_{0} & 0\\
0 & \mu_{1}
\end{array}\right)\right]\cdot\left[\frac{\left(\begin{array}{rr}
u_{11} & -u_{10}\\
-u_{01} & u_{00}
\end{array}\right)}{u_{00}u_{11}-u_{01}u_{10}}\cdot\left(\begin{array}{cc}
\nu_{0} & 0\\
0 & \nu_{1}
\end{array}\right)\right]^{\top}\label{Matrix spectral decomposition}
\end{equation}
The spectral constraints yield generators for the polynomial ideal
$\mathcal{I}_{\mathbf{A}}$ in the polynomial ring\\
$\mathbb{C}\left[u_{00},u_{01},u_{10},u_{11},\frac{u_{00}}{u_{00}u_{11}-u_{01}u_{10}},\frac{u_{01}}{u_{00}u_{11}-u_{01}u_{10}},\frac{u_{10}}{u_{00}u_{11}-u_{01}u_{10}},\frac{u_{11}}{u_{00}u_{11}-u_{01}u_{10}},\mu_{0},\mu_{1},\nu_{0},\nu_{1}\right]$.
The spectral elimination ideal is defined as 
\[
\mathcal{I}_{\mathbf{A}}\cap\mathbb{C}\left[\mu_{0},\mu_{1},\nu_{0},\nu_{1}\right]
\]
The spectral decomposition constraints can thus be rewritten as 
\[
\left(\left(\begin{array}{cccc}
1 & 0 & 0 & 0\\
0 & 1 & 0 & 0\\
0 & 0 & 1 & 0\\
0 & 0 & 0 & 1
\end{array}\right)\otimes\left(\begin{array}{rr}
1 & 1\\
\mu_{0}\nu_{0} & \mu_{1}\nu_{1}
\end{array}\right)\right)\cdot\left(\begin{array}{c}
\frac{u_{00}\cdot u_{11}}{u_{00}u_{11}-u_{01}u_{10}}\\
\frac{-u_{01}\cdot u_{10}}{u_{00}u_{11}-u_{01}u_{10}}\\
\frac{-u_{00}\cdot u_{01}}{u_{00}u_{11}-u_{01}u_{10}}\\
\frac{u_{01}\cdot u_{00}}{u_{00}u_{11}-u_{01}u_{10}}\\
\frac{u_{10}\cdot u_{11}}{u_{00}u_{11}-u_{01}u_{10}}\\
\frac{-u_{11}\cdot u_{10}}{u_{00}u_{11}-u_{01}u_{10}}\\
\frac{-u_{10}\cdot u_{01}}{u_{00}u_{11}-u_{01}u_{10}}\\
\frac{u_{11}\cdot u_{00}}{u_{00}u_{11}-u_{01}u_{10}}
\end{array}\right)=\left(\begin{array}{c}
1\\
a_{00}\\
0\\
a_{01}\\
0\\
a_{10}\\
1\\
a_{11}
\end{array}\right),
\]
from which it follows that 
\[
\left(\begin{array}{c}
\frac{u_{00}\cdot u_{11}}{u_{00}u_{11}-u_{01}u_{10}}\\
\frac{-u_{01}\cdot u_{10}}{u_{00}u_{11}-u_{01}u_{10}}\\
\frac{-u_{00}\cdot u_{01}}{u_{00}u_{11}-u_{01}u_{10}}\\
\frac{u_{01}\cdot u_{00}}{u_{00}u_{11}-u_{01}u_{10}}\\
\frac{u_{10}\cdot u_{11}}{u_{00}u_{11}-u_{01}u_{10}}\\
\frac{-u_{11}\cdot u_{10}}{u_{00}u_{11}-u_{01}u_{10}}\\
\frac{-u_{10}\cdot u_{01}}{u_{00}u_{11}-u_{01}u_{10}}\\
\frac{u_{11}\cdot u_{00}}{u_{00}u_{11}-u_{01}u_{10}}
\end{array}\right)=\left(\left(\begin{array}{cccc}
1 & 0 & 0 & 0\\
0 & 1 & 0 & 0\\
0 & 0 & 1 & 0\\
0 & 0 & 0 & 1
\end{array}\right)\otimes\left(\begin{array}{rr}
1 & 1\\
\mu_{0}\nu_{0} & \mu_{1}\nu_{1}
\end{array}\right)^{-1}\right)\cdot\left(\begin{array}{c}
1\\
a_{00}\\
0\\
a_{01}\\
0\\
a_{10}\\
1\\
a_{11}
\end{array}\right).
\]
Consequently, the entries of the vectors $\left(\begin{array}{c}
\frac{u_{00}u_{11}}{u_{00}u_{11}-u_{01}u_{10}}\\
\frac{-u_{01}u_{10}}{u_{00}u_{11}-u_{01}u_{10}}
\end{array}\right)$, $\left(\begin{array}{c}
\frac{-u_{00}u_{01}}{u_{00}u_{11}-u_{01}u_{10}}\\
\frac{u_{01}u_{00}}{u_{00}u_{11}-u_{01}u_{10}}
\end{array}\right)$, $\left(\begin{array}{c}
\frac{u_{10}u_{11}}{u_{00}u_{11}-u_{01}u_{10}}\\
\frac{-u_{11}u_{10}}{u_{00}u_{11}-u_{01}u_{10}}
\end{array}\right)$, $\left(\begin{array}{c}
\frac{-u_{10}u_{01}}{u_{00}u_{11}-u_{01}u_{10}}\\
\frac{u_{11}u_{00}}{u_{00}u_{11}-u_{01}u_{10}}
\end{array}\right)$\\
\\
can be expressed as rational functions in the variables $\mu_{0}\nu_{0}$
and $\mu_{1}\nu_{1}$. The variables $u_{00}$, $u_{01}$, $u_{10}$,
$u_{11}$ are further eliminated via the algebraic relation 
\[
\left(\begin{array}{c}
\frac{\left(u_{00}u_{11}\right)\left(-u_{10}u_{01}\right)}{\left(u_{00}u_{11}-u_{01}u_{10}\right)^{2}}\\
\frac{\left(-u_{01}u_{10}\right)\left(u_{11}u_{00}\right)}{\left(u_{00}u_{11}-u_{01}u_{10}\right)^{2}}
\end{array}\right)=\left(\begin{array}{c}
\frac{\left(-u_{00}u_{01}\right)\left(u_{10}u_{11}\right)}{\left(u_{00}u_{11}-u_{01}u_{10}\right)^{2}}\\
\frac{\left(u_{01}u_{00}\right)\left(-u_{11}u_{10}\right)}{\left(u_{00}u_{11}-u_{01}u_{10}\right)^{2}}
\end{array}\right).
\]
The algebraic relation above yields the characteristic polynomial
\begin{equation}
\left(\begin{array}{c}
\frac{\left(\mu_{1}\nu_{1}-a_{00}\right)\left(\mu_{1}\nu_{1}-a_{00}\right)}{\left(\mu_{1}\nu_{1}-\mu_{0}\nu_{0}\right)^{2}}\\
\frac{\left(a_{00}-\mu_{0}\nu_{0}\right)\left(a_{11}-\mu_{0}\nu_{0}\right)}{\left(\mu_{1}\nu_{1}-\mu_{0}\nu_{0}\right)^{2}}
\end{array}\right)=\left(\begin{array}{c}
\frac{\left(-a_{01}\right)\left(-a_{10}\right)}{\left(\mu_{1}\nu_{1}-\mu_{0}\nu_{0}\right)^{2}}\\
\frac{a_{01}a_{10}}{\left(\mu_{1}\nu_{1}-\mu_{0}\nu_{0}\right)^{2}}
\end{array}\right).
\end{equation}
Once the determinant polynomial is derived, the generator of the spectral
elimination ideal is more simply obtained by considering the polynomial
\[
\det\left(\mathbf{A}-\mu\nu\,\mathbf{I}_{n}\right).
\]
In particular, in the case $n=2$ we have 
\[
\det\left(\mathbf{A}-\mu\nu\,\mathbf{I}_{n}\right)=\left(\mu\nu\right)^{2}-\mbox{Tr}\left(\mathbf{A}\right)\left(\mu\nu\right)+\det\left(\mathbf{A}\right).
\]
We point out this well-known fact only to emphasize the close analogy
with the hypermatrix case discussed in the next section.\\

\textbf{Theorem 6} : Let $\mathbf{A}\in\mathbb{C}^{n\times n}$ be
a matrix generated by arbitrary combinations of direct sums and Kronecker
products of $2\times2$ matrices. Furthermore, assume that each $2\times2$
generator matrix admits a spectral decomposition. Then $\mathbf{A}$
admits a spectral decomposition of the form 
\[
\mathbf{A}=\left(\mathbf{U}\cdot\mbox{diag}\left(\boldsymbol{\mu}\right)\right)\cdot\left(\left(\mathbf{U}^{-1}\right)^{\top}\cdot\mbox{diag}\left(\boldsymbol{\nu}\right)\right)^{\top}
\]
\emph{Proof}: From the fact that each $2\times2$ generator matrix
admits a spectral decomposition, it follows that the spectral decomposition
of $\mathbf{A}$ is obtained from the spectral decomposition of the
generator matrices by repeated use of Lemma 3.

\subsection{The hypermatrix case}

The spectral decomposition of a hypermatrix $\mathbf{A}\in\mathbb{C}^{2\times2\times2}$
is expressed in terms of an uncorrelated triple ($\mathbf{U}$,$\mathbf{V}$,$\mathbf{W}$).
The $2\times1\times2$ hypermatrix column slices $\left\{ \mathbf{U}\left[:,k,:\right],\,\mathbf{V}\left[:,k,:\right],\:\mathbf{W}\left[:,k,:\right]\right\} _{0\le k<2}$
collect the ``\emph{eigenmatrices}'' of $\mathbf{A}$. We recall
from \cite{GER} that the spectral decomposition is expressed as 
\begin{equation}
\mathbf{A}=\mbox{Prod}\left(\mbox{Prod}\left(\mathbf{U},\mathbf{D}_{0},\mbox{\ensuremath{\mathbf{D}}}_{0}^{\top}\right),\mbox{Prod}\left(\mathbf{V},\mathbf{D}_{1},\mbox{\ensuremath{\mathbf{D}}}_{1}^{\top}\right)^{\top^{2}},\mbox{Prod}\left(\mathbf{W},\mathbf{D}_{2},\mbox{\ensuremath{\mathbf{D}}}_{2}^{\top}\right)^{\top}\right),\label{Hypermatrix spectral decomposition}
\end{equation}
where the $2\times2\times2$ hypermatrices $\mathbf{D}_{0}$, $\mathbf{D}_{1}$,
and $\mathbf{D}_{2}$ are third-order analogs of the diagonal matrices
\[
\left(\begin{array}{cc}
\mu_{0} & 0\\
0 & \mu_{1}
\end{array}\right),\:\left(\begin{array}{cc}
\nu_{0} & 0\\
0 & \nu_{1}
\end{array}\right)
\]
used in \eqref{Matrix spectral decomposition}. The entries of the
hypermatrices $\mathbf{D}_{0}$, $\mathbf{D}_{1}$, and $\mathbf{D}_{2}$
are respectively given by 
\[
\mathbf{D}_{0}\left[:,:,0\right]=\left(\begin{array}{cc}
\mu_{00} & 0\\
\mu_{01} & 0
\end{array}\right),\;\mathbf{D}_{0}\left[:,:,1\right]=\left(\begin{array}{cc}
0 & \mu_{01}\\
0 & \mu_{11}
\end{array}\right),
\]
\[
\mathbf{D}_{1}\left[:,:,0\right]=\left(\begin{array}{cc}
\nu_{00} & 0\\
\nu_{01} & 0
\end{array}\right),\;\mathbf{D}_{1}\left[:,:,1\right]=\left(\begin{array}{cc}
0 & \nu_{01}\\
0 & \nu_{11}
\end{array}\right),
\]
\[
\mathbf{D}_{2}\left[:,:,0\right]=\left(\begin{array}{cc}
\omega_{00} & 0\\
\omega_{01} & 0
\end{array}\right),\;\mathbf{D}_{2}\left[:,:,1\right]=\left(\begin{array}{cc}
0 & \omega_{01}\\
0 & \omega_{11}
\end{array}\right).
\]
The spectral constraints yield generators for the polynomial ideal
$\mathcal{I}_{\mathbf{A}}$ in the polynomial ring\\
$\mathbb{C}\left[u_{000},\cdots,u_{111},\,v_{000},\cdots,v_{111},\,w_{000},\cdots,w_{111},\,\mu_{00},\mu_{01},\mu_{11},\,\nu_{00},\nu_{01},\nu_{11},\,\omega_{00},\omega_{01},\omega_{11}\right]$.\\
\\
By analogy to the matrix derivation, generators for the spectral elimination
ideal are generators for the polynomial ideal 
\[
\mathcal{I}_{\mathbf{A}}\cap\mathbb{C}\left[\mu_{00},\mu_{01},\mu_{11},\,\nu_{00},\nu_{01},\nu_{11},\,\omega_{00},\omega_{01},\omega_{11}\right].
\]
The generators of the elimination ideal suggests the $2\times2\times2$
analog of the determinant as well as the corresponding characteristic
polynomial. We rewrite the hypermatrix spectral decomposition constraints
\eqref{Hypermatrix spectral decomposition} as follows: 
\[
\left[\bigoplus_{0\le i,j,k<2}\left(\mathbf{I}_{2}\otimes\left(\begin{array}{rr}
1 & 1\\
\mu_{0i}\mu_{0k}\nu_{0j}\nu_{0i}\omega_{0k}\omega_{0j} & \mu_{i1}\mu_{k1}\nu_{j1}\nu_{i1}\omega_{k1}\omega_{j1}
\end{array}\right)\right)\right]\left(\begin{array}{c}
u_{00k}\cdot v_{000}\cdot w_{000}\\
u_{010}\cdot v_{010}\cdot w_{010}\\
\vdots\\
u_{i0k}\cdot v_{j0i}\cdot w_{k0j}\\
u_{i1k}\cdot v_{j1i}\cdot w_{k1j}\\
\vdots\\
u_{101}\cdot v_{101}\cdot w_{101}\\
u_{111}\cdot v_{111}\cdot w_{111}
\end{array}\right)=\left(\begin{array}{c}
1\\
a_{000}\\
0\\
a_{001}\\
0\\
a_{010}\\
0\\
a_{011}\\
0\\
a_{100}\\
0\\
a_{101}\\
0\\
a_{110}\\
1\\
a_{111}
\end{array}\right).
\]
It therefore follows from the equality above that 
\[
\left(\begin{array}{c}
u_{00k}\cdot v_{000}\cdot w_{000}\\
u_{010}\cdot v_{010}\cdot w_{010}\\
\vdots\\
u_{i0k}\cdot v_{j0i}\cdot w_{k0j}\\
u_{i1k}\cdot v_{j1i}\cdot w_{k1j}\\
\vdots\\
u_{101}\cdot v_{101}\cdot w_{101}\\
u_{111}\cdot v_{111}\cdot w_{111}
\end{array}\right)=\left[\bigoplus_{0\le i,j,k<2}\left(\mathbf{I}_{2}\otimes\left(\begin{array}{rr}
1 & 1\\
\mu_{0i}\mu_{0k}\nu_{0j}\nu_{0i}\omega_{0k}\omega_{0j} & \mu_{i1}\mu_{k1}\nu_{j1}\nu_{i1}\omega_{k1}\omega_{j1}
\end{array}\right)\right)\right]^{-1}\cdot\left(\begin{array}{c}
1\\
a_{000}\\
0\\
a_{001}\\
0\\
a_{010}\\
0\\
a_{011}\\
0\\
a_{100}\\
0\\
a_{101}\\
0\\
a_{110}\\
1\\
a_{111}
\end{array}\right),
\]
implicitly assuming that 
\[
0\ne\prod_{0\le i,j,k<2}\left(\mu_{0i}\mu_{0k}\nu_{0j}\nu_{0i}\omega_{0k}\omega_{0j}-\mu_{i1}\mu_{k1}\nu_{j1}\nu_{i1}\omega_{k1}\omega_{j1}\right).
\]
Consequently the entries of the vectors $\left\{ \left(\begin{array}{c}
u_{i0k}\cdot v_{j0i}\cdot w_{k0j}\\
u_{i1k}\cdot v_{j1i}\cdot w_{k1j}
\end{array}\right)\right\} _{0\le i,j,k<2}$ are rational functions of the variables $\mu_{00},\mu_{01},\mu_{11},\,\nu_{00},\nu_{01},\nu_{11},\,\omega_{00},\omega_{01},\omega_{11}$.
The variables $u_{000},\cdots,u_{111},\,v_{000},\cdots,v_{111},\,w_{000},\cdots,w_{111}$
are thus eliminated via the relation
\[
\left(\begin{array}{c}
\left(u_{000}v_{000}w_{000}\right)\cdot\left(u_{001}v_{100}w_{101}\right)\cdot\left(u_{101}v_{001}w_{100}\right)\cdot\left(u_{100}v_{101}w_{001}\right)\\
\left(u_{010}v_{010}w_{010}\right)\cdot\left(u_{011}v_{110}w_{111}\right)\cdot\left(u_{111}v_{011}w_{110}\right)\cdot\left(u_{100}v_{101}w_{001}\right)
\end{array}\right)
\]
\[
=
\]
\[
\left(\begin{array}{c}
\left(u_{001}v_{000}w_{100}\right)\cdot\left(u_{000}v_{100}w_{001}\right)\cdot\left(u_{100}v_{001}w_{000}\right)\cdot\left(u_{101}v_{101}w_{101}\right)\\
\left(u_{011}v_{010}w_{110}\right)\cdot\left(u_{010}v_{110}w_{011}\right)\cdot\left(u_{110}v_{011}w_{010}\right)\cdot\left(u_{101}v_{101}w_{101}\right)
\end{array}\right),
\]
which yields the third order analog of the characteristic polynomial
\[
\left(\begin{array}{c}
0\\
0
\end{array}\right)=\prod_{0\le i,j,k<2}\left(\mu_{0i}\mu_{0k}\nu_{0j}\nu_{0i}\omega_{0k}\omega_{0j}-\mu_{i1}\mu_{k1}\nu_{j1}\nu_{i1}\omega_{k1}\omega_{j1}\right)^{-2}\times
\]
\[
\left(\begin{array}{c}
a_{001}a_{010}a_{100}\left(\mu_{11}\nu_{11}\omega_{11}\right)^{2}-a_{011}a_{101}a_{110}\left(\mu_{01}\nu_{01}\omega_{01}\right)^{2}+a_{000}a_{011}a_{101}a_{110}-a_{001}a_{010}a_{100}a_{111}\\
a_{001}a_{010}a_{100}\left(\mu_{01}\nu_{01}\omega_{01}\right)^{2}-a_{011}a_{101}a_{110}\left(\mu_{00}\nu_{00}\omega_{00}\right)^{2}+a_{000}a_{011}a_{101}a_{110}-a_{001}a_{010}a_{100}a_{111}
\end{array}\right).
\]
The generators for the spectral elimination ideal correspond to generators
for the polynomial ideals 
\[
\mathcal{I}_{\mathbf{A}}\cap\mathbb{C}\left[\mu_{00},\mu_{01},\,\nu_{00},\nu_{01},\,\omega_{00},\omega_{01}\right]\:\mbox{ and }\mathcal{I}_{\mathbf{A}}\cap\mathbb{C}\left[\mu_{00},\mu_{01},\,\nu_{00},\nu_{01},\,\omega_{00},\omega_{01}\right]
\]
respectively given by 
\[
a_{001}a_{010}a_{100}\left(\mu_{01}\nu_{01}\omega_{01}\right)^{2}-a_{011}a_{101}a_{110}\left(\mu_{00}\nu_{00}\omega_{00}\right)^{2}+a_{000}a_{011}a_{101}a_{110}-a_{001}a_{010}a_{100}a_{111}
\]
\[
\mbox{and}
\]
\[
a_{001}a_{010}a_{100}\left(\mu_{11}\nu_{11}\omega_{11}\right)^{2}-a_{011}a_{101}a_{110}\left(\mu_{01}\nu_{01}\omega_{01}\right)^{2}+a_{000}a_{011}a_{101}a_{110}-a_{001}a_{010}a_{100}a_{111}.
\]
\\
The derivation also suggests that the $2\times2\times2$ hypermatrix
analog of the characteristic polynomial is the polynomial
\[
p\left(\mu_{0}\nu_{0}\omega_{0},\:\mu_{1}\nu_{1}\omega_{1}\right)=
\]
\[
a_{001}a_{010}a_{100}\left(\mu_{1}\nu_{1}\omega_{1}\right)^{2}-a_{011}a_{101}a_{110}\left(\mu_{0}\nu_{0}\omega_{0}\right)^{2}+\left(a_{000}a_{011}a_{101}a_{110}-a_{001}a_{010}a_{100}a_{111}\right)
\]
whose constant term $a_{000}a_{011}a_{101}a_{110}-a_{001}a_{010}a_{100}a_{111}$
is the $2\times2\times2$ hypermatrix analog of the \emph{determinant}
polynomial. Note that the determinant polynomial is linear in the
row, column and depth slices. Furthermore, the $2\times2\times2$
analog of the determinant changes sign with a row, column or depth
slice exchange. Finally, the determinant of a $2\times2\times2$ hypermatrix
non zero if and only if the BM-rank of corresponding hypermatrix is
equal to 2.\textbf{}\\
\textbf{}\\
\textbf{}\\
\textbf{Theorem 7} : Let $\mathbf{A}\in\mathbb{C}^{n\times n\times n}$
be a hypermatrix generated by some arbitrary combination of direct
sums and Kronecker products of $2\times2\times2$ hypermatrices. Furthermore,
assume that each $2\times2\times2$ generator hypermatrix admits a
spectral decomposition. Then $\mathbf{A}$ admits a spectral decomposition
of the form 
\[
\mathbf{A}=\mbox{Prod}\left(\mbox{Prod}\left(\mathbf{U},\mathbf{D}_{0},\mbox{\ensuremath{\mathbf{D}}}_{0}^{\top}\right),\mbox{Prod}\left(\mathbf{V},\mathbf{D}_{1},\mbox{\ensuremath{\mathbf{D}}}_{1}^{\top}\right)^{\top^{2}},\mbox{Prod}\left(\mathbf{W},\mathbf{D}_{2},\mbox{\ensuremath{\mathbf{D}}}_{2}^{\top}\right)^{\top}\right),
\]
\[
\mbox{subject to }
\]
\[
\mbox{Prod}\left(\mathbf{U},\mathbf{V}^{\top^{2}},\mathbf{W}^{\top}\right)=\boldsymbol{\Delta}
\]
\[
\left[\mathbf{D}_{0}\right]_{ijk}=\begin{cases}
\begin{array}{cc}
\mu_{jk}=\mu_{kj} & \mbox{if }0\le i=k<n\\
0 & \mbox{otherwise }
\end{array}\end{cases},
\]
\[
\left[\mathbf{D}_{1}\right]_{ijk}=\begin{cases}
\begin{array}{cc}
\nu_{jk}=\nu_{kj} & \mbox{if }0\le i=k<n\\
0 & \mbox{otherwise }
\end{array}\end{cases},
\]
\[
\left[\mathbf{D}_{2}\right]_{ijk}=\begin{cases}
\begin{array}{cc}
\omega_{jk}=\omega_{kj} & \mbox{if }0\le i=k<n\\
0 & \mbox{otherwise }
\end{array}\end{cases}.
\]
\emph{Proof}: From the fact that each $2\times2\times2$ generator
hypermatrix admits a spectral decomposition, It follows that the spectral
decomposition of $\mathbf{A}$ is derived from the spectral decomposition
of the generators by repeated use of Lemma 3.\\
\\
\\
As in the matrix case, the characteristic polynomial can be obtained
directly from the $2\times2\times2$ analog of the determinant polynomial
derived above for the cubic side length $2$ hypermatrix $\mathbf{B}$
whose entries are given by 
\[
\left[\mathbf{B}\right]_{ijk}=a_{ijk}-\sum_{0\le t<2}\left(\mu_{i}u_{itk}\mu_{k}\right)\left(\nu_{j}v_{jti}\nu_{i}\right)\left(\omega_{k}w_{ktj}\omega_{j}\right)
\]
\[
\Rightarrow\left[\mathbf{B}\right]_{ijk}=a_{ijk}-\mu_{i}\mu_{k}\nu_{j}\nu_{i}\omega_{k}\omega_{j}\sum_{0\le t<2}u_{itk}v_{jti}w_{ktj}.
\]
From the fact that $\mbox{Prod}\left(\mathbf{U},\mathbf{V},\mathbf{W}\right)=\boldsymbol{\Delta}$
it follows that 
\[
\left[\mathbf{B}\right]_{ijk}=\begin{cases}
\begin{array}{cc}
a_{iii}-\left(\mu_{i}\nu_{i}\omega_{i}\right)^{2} & \mbox{ if }\:0\le i=j=k<2\\
a_{ijk} & \mbox{otherwise}
\end{array}\end{cases}.
\]
The characteristic polynomial is thus obtained by computing the $2\times2\times2$
analog of the determinant associated with $\mathbf{B}$ 
\[
\mbox{det}\left(\mathbf{B}\right)=
\]
\begin{equation}
a_{001}a_{010}a_{100}\left(\mu_{1}\nu_{1}\omega_{1}\right)^{2}-a_{011}a_{101}a_{110}\left(\mu_{0}\nu_{0}\omega_{0}\right)^{2}+\left(a_{000}a_{011}a_{101}a_{110}-a_{001}a_{010}a_{100}a_{111}\right).
\end{equation}
\\
\\
The $m$-th order side length $2$ analog of the determinant is derived
in a similar way from the family of spectral elimination ideals
\[
\mbox{det}\left(\mathbf{A}\right)=
\]
\begin{equation}
\left(\prod_{\begin{array}{c}
\mathbf{j}\in\left\{ 0,1\right\} ^{1\times m}\\
\left\Vert \mathbf{j}\right\Vert _{\ell_{1}}\equiv0\mod2
\end{array}}a_{\mathbf{j}}\right)-\left(\prod_{\begin{array}{c}
\mathbf{j}\in\left\{ 0,1\right\} ^{1\times m}\\
\left\Vert \mathbf{j}\right\Vert _{\ell_{1}}\equiv1\mod2
\end{array}}a_{\mathbf{j}}\right),\label{Hyperdeterminant}
\end{equation}
for any order $m$ hypermatrix $\mathbf{A}$ with side length $2$.\\
\\
We remark that the spectral decomposition described here is different
from the approaches first introduced by Liqun Qi and Lek-Heng Lim
in \cite{Lim,Qi:2005:ERS:1740736.1740799}. The first essential distinction
arises from the fact that their proposed generalization to hypermatrices/tensors
of the notion of eigenvalues is not associated with any particular
hypermatrix factorization, although it suggests various rank one approximation
schemes. The second distinction arises from the fact that the E-characteristic
polynomial is defined for hypermatrices which are symmetric relative
to any permutation of the entries, whereas our proposed formulation
makes no such restrictions.

\subsection{Spectra of adjacency hypermatrices of groups}

As an illustration of naturally occurring hypermatrices we consider
the adjacency hypermatrices of finite groups. To an arbitrary finite
group $G$ of order $n$, one associates an $n\times n\times n$ adjacency
hypermatrix $\mathbf{A}_{G}$ with binary entries specified as follows:
\[
\forall\:i,j,k\in G,\quad a_{ijk}=\begin{cases}
\begin{array}{cc}
1 & \mbox{if }i\cdot j=k\:\mbox{ in }G\\
0 & \mbox{otherwise}
\end{array}\end{cases}.
\]
As illustration, we consider here adjacency hypermatrices associated
with the family of groups of the form $\nicefrac{\mathbb{Z}}{2\mathbb{Z}}\times\nicefrac{\mathbb{Z}}{2\mathbb{Z}}\times\cdots\times\nicefrac{\mathbb{Z}}{2\mathbb{Z}}$.
Note that by definition 
\[
\mathbf{A}_{\nicefrac{\mathbb{Z}}{2\mathbb{Z}}\times\nicefrac{\mathbb{Z}}{2\mathbb{Z}}\times\cdots\times\nicefrac{\mathbb{Z}}{2\mathbb{Z}}}=\mathbf{A}_{\nicefrac{\mathbb{Z}}{2\mathbb{Z}}}\otimes\mathbf{A}_{\nicefrac{\mathbb{Z}}{2\mathbb{Z}}}\otimes\cdots\otimes\mathbf{A}_{\nicefrac{\mathbb{Z}}{2\mathbb{Z}}}.
\]
Consequently, the spectral decomposition of the adjacency hypermatrix
$\mathbf{A}_{\nicefrac{\mathbb{Z}}{2\mathbb{Z}}\times\nicefrac{\mathbb{Z}}{2\mathbb{Z}}\times\cdots\times\nicefrac{\mathbb{Z}}{2\mathbb{Z}}}$
is determined by the spectral decomposition of the $2\times2\times2$
hypermatrix $\mathbf{A}_{\nicefrac{\mathbb{Z}}{2\mathbb{Z}}}$. The
entries of the hypermatrix $\mathbf{A}_{\nicefrac{\mathbb{Z}}{2\mathbb{Z}}}$
are given by
\[
\forall\:i,j,k\in\nicefrac{\mathbb{Z}}{2\mathbb{Z}},\quad\left[\mathbf{A}_{\nicefrac{\mathbb{Z}}{2\mathbb{Z}}}\right]_{i,j,k}=\begin{cases}
\begin{array}{cc}
1 & \mbox{if }i+j\equiv k\mod2\\
0 & \mbox{otherwise}
\end{array}\end{cases},
\]
\[
\mathbf{A}_{\nicefrac{\mathbb{Z}}{2\mathbb{Z}}}\left[:,:,0\right]=\left(\begin{array}{cc}
1 & 0\\
0 & 1
\end{array}\right),\quad\mathbf{A}_{\nicefrac{\mathbb{Z}}{2\mathbb{Z}}}\left[:,:,1\right]=\left(\begin{array}{cc}
0 & 1\\
1 & 0
\end{array}\right).
\]
By symmetry the hypermatrix $\mathbf{A}$ admits a spectral decomposition
of the form 
\[
\mathbf{A}_{\nicefrac{\mathbb{Z}}{2\mathbb{Z}}}=\mbox{Prod}\left(\mbox{Prod}\left(\mathbf{Q},\mathbf{D},\mathbf{D}^{\top}\right),\mbox{Prod}\left(\mathbf{Q},\mathbf{D},\mathbf{D}^{\top}\right)^{\top^{2}},\mbox{Prod}\left(\mathbf{Q},\mathbf{D},\mathbf{D}^{\top}\right)^{\top}\right),
\]
where the hypermatrix $\mathbf{D}$ is of the form 
\[
\mathbf{D}_{0}\left[:,:,0\right]=\left(\begin{array}{cc}
\lambda_{00} & 0\\
\lambda_{01} & 0
\end{array}\right),\;\mathbf{D}_{0}\left[:,:,1\right]=\left(\begin{array}{cc}
0 & \lambda_{01}\\
0 & \lambda_{11}
\end{array}\right),
\]
and the hypermatrix $\mathbf{Q}$ is subject to the orthogonality
constraints expressed by 
\[
\mbox{Prod}\left(\mathbf{Q},\mathbf{Q}^{\top^{2}},\mathbf{Q}^{\top}\right)=\boldsymbol{\Delta}.
\]
The spectrum of $\mathbf{A}_{\nicefrac{\mathbb{Z}}{2\mathbb{Z}}}$
is determined by the following parametrization of orthogonal hypermatrices
\[
\mathbf{Q}\left[:,:,0\right]=\left(\begin{array}{rr}
\left(x^{3}+1\right)^{-\frac{1}{3}} & \left(\frac{1}{x^{3}}+1\right)^{-\frac{1}{3}}\\
-x & 1
\end{array}\right),\mathbf{Q}\left[:,:,1\right]=\left(\begin{array}{rr}
1 & 1\\
\left(x^{3}+1\right)^{-\frac{1}{3}} & \left(\frac{1}{x^{3}}+1\right)^{-\frac{1}{3}}
\end{array}\right),
\]
as well as the following parametrization for the hypermatrix $\mathbf{D}$:
\[
\mathbf{D}\left[:,:,0\right]=\left(\begin{array}{rr}
\left(-x^{3}\right)^{\frac{1}{12}} & 0\\
\left(-x^{3}\right)^{\frac{1}{6}} & 0
\end{array}\right),\:\mathbf{D}\left[:,:,1\right]=\left(\begin{array}{rr}
0 & \left(-x^{3}\right)^{\frac{1}{6}}\\
0 & 1
\end{array}\right).
\]
The parametrization above found via \cite{sage} ensures that 
\[
\forall\:\left(i,j,k\right)\in\left\{ \left(1,0,0\right),\left(0,1,0\right),\left(0,0,1\right),\left(1,1,1\right)\right\} ,
\]
\[
\left[\mbox{Prod}\left(\mbox{Prod}\left(\mathbf{Q},\mathbf{D},\mathbf{D}^{\top}\right),\mbox{Prod}\left(\mathbf{Q},\mathbf{D},\mathbf{D}^{\top}\right)^{\top^{2}},\mbox{Prod}\left(\mathbf{Q},\mathbf{D},\mathbf{D}^{\top}\right)^{\top}\right)\right]_{i,j,k}=0.
\]
Finally, by symmetry, the spectral decomposition of $\mathbf{A}$
is obtained by solving for the parameter $x$ in the equation 
\[
\left[\mbox{Prod}\left(\mbox{Prod}\left(\mathbf{Q},\mathbf{D},\mathbf{D}^{\top}\right),\mbox{Prod}\left(\mathbf{Q},\mathbf{D},\mathbf{D}^{\top}\right)^{\top^{2}},\mbox{Prod}\left(\mathbf{Q},\mathbf{D},\mathbf{D}^{\top}\right)^{\top}\right)\right]_{0,0,0}=
\]
\[
\left[\mbox{Prod}\left(\mbox{Prod}\left(\mathbf{Q},\mathbf{D},\mathbf{D}^{\top}\right),\mbox{Prod}\left(\mathbf{Q},\mathbf{D},\mathbf{D}^{\top}\right)^{\top^{2}},\mbox{Prod}\left(\mathbf{Q},\mathbf{D},\mathbf{D}^{\top}\right)^{\top}\right)\right]_{0,1,1},
\]
which yields the equation 
\[
\frac{\left(-1\right)^{\frac{5}{6}}x^{\frac{7}{2}}-\left(-1\right)^{\frac{1}{3}}x^{2}}{\left(x^{2}-x+1\right)^{\frac{1}{3}}\left(x+1\right)^{\frac{1}{3}}}-\frac{x^{6}+x\sqrt{-x}}{x^{3}+1}=0,
\]
for which the existence of complex roots follows immediately from
the fundamental theorem of algebra. Consequently, by Lemma 3 the spectral
decomposition of the $m$-th order adjacency hypermatrix of the group
$\nicefrac{\mathbb{Z}}{2\mathbb{Z}}\times\nicefrac{\mathbb{Z}}{2\mathbb{Z}}\times\cdots\times\nicefrac{\mathbb{Z}}{2\mathbb{Z}}$
is expressed as
\[
\mathbf{A}_{\nicefrac{\mathbb{Z}}{2\mathbb{Z}}\times\nicefrac{\mathbb{Z}}{2\mathbb{Z}}\times\cdots\times\nicefrac{\mathbb{Z}}{2\mathbb{Z}}}=\left(\mathbf{A}_{\nicefrac{\mathbb{Z}}{2\mathbb{Z}}}\right)^{\otimes^{n}}=
\]
\[
\mbox{Prod}\left(\mbox{Prod}\left(\mathbf{Q}^{\otimes^{n}},\mathbf{D}^{\otimes^{n}},\left(\mathbf{D}^{\otimes^{n}}\right)^{\top}\right),\mbox{Prod}\left(\mathbf{Q}^{\otimes^{n}},\mathbf{D}^{\otimes^{n}},\left(\mathbf{D}^{\otimes^{n}}\right)^{\top}\right)^{\top^{2}},\mbox{Prod}\left(\mathbf{Q}^{\otimes^{n}},\mathbf{D}^{\otimes^{n}},\left(\mathbf{D}^{\otimes^{n}}\right)^{\top}\right)^{\top}\right).
\]

\section{General matrix and hypermatrix Rayleigh quotient}

The Rayleigh quotient is central to many applications of the spectral
decomposition of matrices. We prove here a slight generalization of
the matrix Rayleigh quotient inequalities. The proposed variant of
the Rayleigh quotient inequalities does not assume Hermicity of the
underlying matrix. We also extend the result to hypermatrices.\\
\textbf{}\\
\textbf{Theorem 8} : Let $\mathbf{A}\in\mathbb{C}^{n\times n}$ having
non-negative eigenvalues. Let the spectral decomposition of $\mathbf{A}$
be given by 
\[
\mathbf{A}=\mathbf{U}\cdot\mbox{diag}\left(\begin{array}{c}
\lambda_{0}\\
\vdots\\
\lambda_{n-1}
\end{array}\right)\cdot\mathbf{V}^{\top},\,\mbox{subject to}\quad\mathbf{I}_{n}=\mathbf{U}\cdot\mathbf{V}^{\top},
\]
Let $\mathbf{P}_{k}=\mbox{Prod}_{\boldsymbol{\Delta}^{(k)}}\left(\mathbf{U},\mathbf{V}^{\top}\right)$
and $S_{k}\subset\mathbb{C}^{n\times1}\times\mathbb{C}^{n\times1}$
be such that $\forall\:\left(\mathbf{x},\mathbf{y}\right)\in S_{k}$,
$\mbox{Prod}_{\mathbf{P}_{k}}\left(\mathbf{x}^{\top},\mathbf{y}\right)\ge0$,
then 
\[
\forall\,\left(\mathbf{x},\mathbf{y}\right)\in\bigcap_{0\le k<n}S_{k},\quad\min_{0\le t<n}\lambda_{t}\le\frac{\mbox{Prod}_{\mathbf{A}}\left(\mathbf{x}^{\top},\mathbf{y}\right)}{\mbox{Prod}\left(\mathbf{x}^{\top},\mathbf{y}\right)}\le\max_{0\le t<n}\lambda_{t}.
\]
( Assuming that $\mbox{Prod}\left(\mathbf{x}^{\top},\mathbf{y}\right)\ne0$
)\\
\\
\emph{Proof} : From the spectral decomposition of $\mathbf{A}$ we
have 
\[
\mbox{Prod}_{\mathbf{A}}\left(\mathbf{x}^{\top},\mathbf{y}\right)=\sum_{0\le k<n}\lambda_{k}\,\mbox{Prod}_{\mathbf{P}_{k}}\left(\mathbf{x}^{\top},\mathbf{y}\right).
\]
By positivity we have 
\[
\forall\,\left(\mathbf{x},\mathbf{y}\right)\in\bigcap_{0\le k<n}S_{k},
\]
\[
\sum_{0\le k<n}\left(\min_{0\le t<n}\lambda_{t}\right)\ \mbox{Prod}_{\mathbf{P}_{k}}\left(\mathbf{x}^{\top},\mathbf{y}\right)\le\mbox{Prod}_{\mathbf{A}}\left(\mathbf{x}^{\top},\mathbf{y}\right)\le\sum_{0\le k<n}\left(\max_{0\le t<n}\lambda_{t}\right)\ \mbox{Prod}_{\mathbf{P}_{k}}\left(\mathbf{x}^{\top},\mathbf{y}\right),
\]
\[
\Rightarrow\min_{0\le t<n}\lambda_{t}\:\sum_{0\le k<n}\mbox{Prod}_{\mathbf{P}_{k}}\left(\mathbf{x}^{\top},\mathbf{y}\right)\le\mbox{Prod}_{\mathbf{A}}\left(\mathbf{x}^{\top},\mathbf{y}\right)\le\max_{0\le t<n}\lambda_{t}\:\sum_{0\le k<n}\mbox{Prod}_{\mathbf{P}_{k}}\left(\mathbf{x}^{\top},\mathbf{y}\right),
\]
which follows from the fact that $\forall\,0\le k<n,\;\mbox{Prod}_{\mathbf{P}_{k}}\left(\mathbf{x}^{\top},\mathbf{y}\right)\ge0$.
By the Parseval identity 
\[
\mbox{Prod}\left(\mathbf{x}^{\top},\mathbf{y}\right)=\sum_{0\le k<n}\mbox{Prod}_{\mathbf{P}_{k}}\left(\mathbf{x}^{\top},\mathbf{y}\right),
\]
\[
\Rightarrow\left(\min_{0\le t<n}\lambda_{t}\right)\,\mbox{Prod}\left(\mathbf{x}^{\top},\mathbf{y}\right)\le\mbox{Prod}_{\mathbf{A}}\left(\mathbf{x}^{\top},\mathbf{y}\right)\le\left(\max_{0\le t<n}\lambda_{t}\right)\,\mbox{Prod}\left(\mathbf{x}^{\top},\mathbf{y}\right),
\]
and the sought after result follows 
\[
\min_{0\le t<n}\lambda_{t}\le\frac{\mbox{Prod}_{\mathbf{A}}\left(\mathbf{x}^{\top},\mathbf{y}\right)}{\mbox{Prod}\left(\mathbf{x}^{\top},\mathbf{y}\right)}\le\max_{0\le t<n}\lambda_{t}\,.\square
\]
By sorting the eigenvalues such that $\lambda_{0}\le\lambda_{1}\le\cdots\le\lambda_{n-2}\le\lambda_{n-1}$
It is then easily verified that the bounds are attained for the choices
\[
\mathbf{x}=\mathbf{V}\left[:,\,0\right],\quad\mathbf{y}=\mathbf{U}\left[:,\,0\right]
\]
\[
\mbox{and}
\]
\[
\mathbf{x}=\mathbf{V}\left[:,\,n-1\right],\quad\mathbf{y}=\mathbf{U}\left[:,\,n-1\right].
\]
It is useful to provide some explicit description for vectors in the
set 
\[
\bigcap_{0\le k<n}S_{k}\subset\mathbb{C}^{n\times1}\times\mathbb{C}^{n\times1}.
\]
The explicit description is given by 
\[
\left(\mathbf{x},\mathbf{y}\right)\in\bigcap_{0\le k<n}S_{k}\Leftrightarrow\mathbf{x}=\sum_{0\le i<n}\alpha_{i}\,\mathbf{V}\left[:,i\right],\;\mathbf{y}=\sum_{0\le j<n}\beta_{j}\,\mathbf{U}\left[:,j\right]\;\mbox{s.t.}\;\left\{ \alpha_{k}\,\beta_{k}\right\} _{0\le k<n}\subset\mathbb{R}_{\ge0}.
\]
\\
Having discussed the matrix formulation of the general Rayleigh quotient,
we now discuss the hypermatrix formulation of the Rayleigh quotient.
For notational convenience we restrict the discussion to third order
hypermatrices, but the formulation extends to hypermatrices of all
orders.\\
\\
\textbf{Theorem 9} : Let $\mathbf{A}\in\mathbb{C}^{n\times n\times n}$,
whose spectral decomposition is given by 
\[
\mathbf{A}=\mbox{Prod}\left(\mbox{Prod}\left(\mathbf{U},\mathbf{D}_{0},\mbox{\ensuremath{\mathbf{D}}}_{0}^{\top}\right),\mbox{Prod}\left(\mathbf{V},\mathbf{D}_{1},\mbox{\ensuremath{\mathbf{D}}}_{1}^{\top}\right)^{\top^{2}},\mbox{Prod}\left(\mathbf{W},\mathbf{D}_{2},\mbox{\ensuremath{\mathbf{D}}}_{2}^{\top}\right)^{\top}\right)
\]
\[
\mbox{subject to}
\]
\[
\mathbf{U},\mathbf{V},\mathbf{W}\in\mathbb{C}^{n\times n\times n}\,\mbox{ and }\,\left[\mbox{Prod}\left(\mathbf{U},\mathbf{V}^{\top^{2}},\mathbf{W}^{\top}\right)\right]_{i,j,k}=\begin{cases}
\begin{array}{cc}
1 & \mbox{ if }0\le i=j=k<n\\
0 & \mbox{otherwise}
\end{array}\end{cases},
\]
where 
\[
\left[\mathbf{D}_{0}\right]_{ijk}=\begin{cases}
\begin{array}{cc}
\mu_{jk}=\mu_{kj}\ge0 & \mbox{if }0\le i=k<n\\
0 & \mbox{otherwise }
\end{array}\end{cases},
\]
\[
\left[\mathbf{D}_{1}\right]_{ijk}=\begin{cases}
\begin{array}{cc}
\nu_{jk}=\nu_{kj}\ge0 & \mbox{if }0\le i=k<n\\
0 & \mbox{otherwise }
\end{array}\end{cases},
\]
\[
\left[\mathbf{D}_{2}\right]_{ijk}=\begin{cases}
\begin{array}{cc}
\omega_{jk}=\omega_{kj}\ge0 & \mbox{if }0\le i=k<n\\
0 & \mbox{otherwise }
\end{array}\end{cases}.
\]
Let $\mathbf{P}_{k}=\mbox{Prod}_{\boldsymbol{\Delta}^{(k)}}\left(\mathbf{U},\mathbf{V}^{\top^{2}},\mathbf{W}^{\top}\right)$
and $S_{k}\subset\mathbb{C}^{n\times1\times1}\times\mathbb{C}^{n\times1\times1}\times\mathbb{C}^{n\times1\times1}$
be such that
\[
\forall\:\left(\mathbf{x},\mathbf{y},\mathbf{z}\right)\in S_{k},\quad\mbox{Prod}_{\mathbf{P}_{k}}\left(\mathbf{x}^{\top^{2}},\mathbf{y}^{\top},\mathbf{z}\right)\ge0.
\]
Then $\forall\,\left(\mathbf{x},\mathbf{y},\mathbf{z}\right)\in\bigcap_{0\le k<n}S_{k},$
\[
\min_{0\le i,j,k,t<n}\left(\mu_{it}\mu_{tk}\,\nu_{jt}\nu_{ti}\,\omega_{kt}\omega_{tj}\right)\le\frac{\mbox{Prod}_{\mathbf{A}}\left(\mathbf{x}^{\top^{2}},\mathbf{y}^{\top},\mathbf{z}\right)}{\mbox{Prod}\left(\mathbf{x}^{\top^{2}},\mathbf{y}^{\top},\mathbf{z}\right)}\le\max_{0\le i,j,k,t<n}\left(\omega_{it}\omega_{kt}\,\nu_{jt}\nu_{it}\,\mu_{kt}\mu_{jt}\right)
\]
( Assuming that $\mbox{Prod}\left(\mathbf{x}^{\top^{2}},\mathbf{y}^{\top},\mathbf{z}\right)\ne0$
)\\
\\
\\
\emph{Proof} : The argument is similar to the matrix case. Recall
from the general Parseval identity that 
\[
\forall\,\left(\mathbf{x},\mathbf{y},\mathbf{z}\right)\in\:\mathbb{C}^{n\times1\times1}\times\mathbb{C}^{n\times1\times1}\times\mathbb{C}^{n\times1\times1},
\]
\[
\mbox{Prod}\left(\mathbf{x}^{\top^{2}},\mathbf{y}^{\top},\mathbf{z}\right)=\sum_{0\le k<n}\mbox{Prod}_{\mathbf{P}_{k}}\left(\mathbf{x}^{\top^{2}},\mathbf{y}^{\top},\mathbf{z}\right),
\]
and by positivity we have 
\[
\forall\,\left(\mathbf{x},\mathbf{y},\mathbf{z}\right)\in\bigcap_{0\le k<n}S_{k},
\]
\[
\sum_{0\le t<n}\min_{0\le i,j,k<n}\left(\mu_{it}\mu_{tk}\,\nu_{jt}\nu_{ti}\,\omega_{kt}\omega_{tj}\right)\,\mbox{Prod}_{\mathbf{P}_{t}}\left(\mathbf{x}^{\top^{2}},\mathbf{y}^{\top},\mathbf{z}\right)\le\mbox{Prod}_{\mathbf{A}}\left(\mathbf{x}^{\top^{2}},\mathbf{y}^{\top},\mathbf{z}\right)
\]
\[
\mbox{and}
\]
\[
\mbox{Prod}_{\mathbf{A}}\left(\mathbf{x}^{\top^{2}},\mathbf{y}^{\top},\mathbf{z}\right)\le\sum_{0\le t<n}\max_{0\le i,j,k<n}\left(\mu_{it}\mu_{tk}\,\nu_{jt}\nu_{ti}\,\omega_{kt}\omega_{tj}\right)\,\mbox{Prod}_{\mathbf{P}_{t}}\left(\mathbf{x}^{\top^{2}},\mathbf{y}^{\top},\mathbf{z}\right)
\]
hence 
\[
\min_{0\le i,j,k,t<n}\left(\mu_{it}\mu_{tk}\,\nu_{jt}\nu_{ti}\,\omega_{kt}\omega_{tj}\right)\,\sum_{0\le t<n}\mbox{Prod}_{\mathbf{P}_{t}}\left(\mathbf{x}^{\top^{2}},\mathbf{y}^{\top},\mathbf{z}\right)\le\mbox{Prod}_{\mathbf{A}}\left(\mathbf{x}^{\top^{2}},\mathbf{y}^{\top},\mathbf{z}\right)
\]
\[
\mbox{and}
\]
\[
\mbox{Prod}_{\mathbf{A}}\left(\mathbf{x}^{\top^{2}},\mathbf{y}^{\top},\mathbf{z}\right)\le\max_{0\le i,j,k,t<n}\left(\mu_{it}\mu_{tk}\,\nu_{jt}\nu_{ti}\,\omega_{kt}\omega_{tj}\right)\,\sum_{0\le t<n}\mbox{Prod}_{\mathbf{P}_{t}}\left(\mathbf{x}^{\top^{2}},\mathbf{y}^{\top},\mathbf{z}\right),
\]
which follows from the fact that $\forall\:\left(\mathbf{x},\mathbf{y},\mathbf{z}\right)\in S_{k}$,
$\mbox{Prod}_{\mathbf{P}_{k}}\left(\mathbf{x}^{\top^{2}},\mathbf{y}^{\top},\mathbf{z}\right)\ge0$.
By the Parseval identity 
\[
\mbox{Prod}\left(\mathbf{x}^{\top^{2}},\mathbf{y}^{\top},\mathbf{z}\right)=\sum_{0\le k<n}\mbox{Prod}_{\mathbf{P}_{k}}\left(\mathbf{x}^{\top^{2}},\mathbf{y}^{\top},\mathbf{z}\right),
\]
we have 
\[
\min_{0\le i,j,k,t<n}\left(\mu_{it}\mu_{tk}\,\nu_{jt}\nu_{ti}\,\omega_{kt}\omega_{tj}\right)\,\mbox{Prod}\left(\mathbf{x}^{\top^{2}},\mathbf{y}^{\top},\mathbf{z}\right)\le\mbox{Prod}_{\mathbf{A}}\left(\mathbf{x}^{\top^{2}},\mathbf{y}^{\top},\mathbf{z}\right)
\]
\[
\mbox{and}
\]
\[
\mbox{Prod}_{\mathbf{A}}\left(\mathbf{x}^{\top^{2}},\mathbf{y}^{\top},\mathbf{z}\right)\le\max_{0\le i,j,k,t<n}\left(\mu_{it}\mu_{tk}\,\nu_{jt}\nu_{ti}\,\omega_{kt}\omega_{tj}\right)\,\mbox{Prod}\left(\mathbf{x}^{\top^{2}},\mathbf{y}^{\top},\mathbf{z}\right).
\]
from which we obtain the sought after result 
\[
\min_{0\le i,j,k,t<n}\left(\mu_{it}\mu_{tk}\,\nu_{jt}\nu_{ti}\,\omega_{kt}\omega_{tj}\right)\le\frac{\mbox{Prod}_{\mathbf{A}}\left(\mathbf{x}^{\top^{2}},\mathbf{y}^{\top},\mathbf{z}\right)}{\mbox{Prod}\left(\mathbf{x}^{\top^{2}},\mathbf{y}^{\top},\mathbf{z}\right)}\le\max_{0\le i,j,k,t<n}\left(\omega_{it}\omega_{kt}\,\nu_{jt}\nu_{it}\,\mu_{kt}\mu_{jt}\right).
\]
\\
\\
For practical uses of the hypermatrix formulation of the Rayleigh
quotient it is useful to provide some explicit description for vectors
in the set 
\[
\bigcap_{0\le k<n}S_{k}\subset\mathbb{C}^{n\times1\times1}\times\mathbb{C}^{n\times1\times1}\times\mathbb{C}^{n\times1\times1}.
\]
We provide here such a characterization. Let $\mathbf{U}$, $\mathbf{V}$,
$\mathbf{W}$ denote uncorrelated third order hypermatrix triple.
Let $\mathbf{P}_{k}$ denote the outer product 
\[
\mathbf{P}_{k}=\mbox{Prod}_{\boldsymbol{\Delta}^{(k)}}\left(\mathbf{U},\mathbf{V}^{\top^{2}},\mathbf{W}^{\top}\right).
\]
We first observe that for each $\mathbf{P}_{k}$ there is a unique
matrix $\mathbf{M}_{k}\left(\mathbf{z}\right)$ for which the following
equality holds 
\[
\forall\:0\le k<2,\quad\mbox{Prod}_{\mathbf{P}_{k}}\left(\mathbf{x}^{\top^{2}},\mathbf{y}^{\top},\mathbf{z}\right)=\mathbf{x}^{T}\cdot\mathbf{M}_{k}\left(\mathbf{z}\right)\cdot\mathbf{y}
\]
consequently the vector $\mathbf{z}$ must be chosen if at all possible
to ensure that the $n\times n$ matrix $\mathbf{M}_{k}\left(\mathbf{z}\right)$
is diagonalizable with positive eigenvalues for all $0\le k<n$. In
the special case of direct sum and Kronecker product constructions
generated by side length $2$ hypermatrices the analysis reduces by
Lemma 3 to the case $n=2$ in which case both of these requirement
are met when
\[
\forall\:0\le k<2,
\]
\[
\mbox{Tr}\left(\mathbf{M}_{k}\left(\mathbf{z}\right)\right)^{2}-4\det\left(\mathbf{M}_{k}\left(\mathbf{z}\right)\right)>0,\;\mbox{Tr}\left(\mathbf{M}_{k}\left(\mathbf{z}\right)\right)\ge0\;\mbox{ and }\;\det\left(\mathbf{M}_{k}\left(\mathbf{z}\right)\right)\ge0.
\]
\[
\mbox{or}
\]
\[
\mbox{Tr}\left(\mathbf{M}_{k}\left(\mathbf{z}\right)\right)^{2}-4\det\left(\mathbf{M}_{k}\left(\mathbf{z}\right)\right)=0,\;\mbox{Tr}\left(\mathbf{M}_{k}\left(\mathbf{z}\right)\right)\ge0\;\mbox{ and }\;\mathbf{M}_{k}\left(\mathbf{z}\right)=\mathbf{M}_{k}^{T}\left(\mathbf{z}\right).
\]
Finally, provided that $\mathbf{z}$ is chosen such that $\mathbf{M}_{k}\left(\mathbf{z}\right)$
is diagonalizable with positive eigenvalues for all $0\le k<n$, Theorem
8 provides a complete characterization for the possible vectors $\mathbf{x}$
and $\mathbf{y}$ for each of the sets $S_{k}$.\\
Furthermore, for a symmetric $n\times n\times n$ hypermatrix $\mathbf{A}$
whose spectral decomposition is expressed by 
\[
\mathbf{A}=\mbox{Prod}\left(\mbox{Prod}\left(\mathbf{Q},\mathbf{D},\mathbf{D}^{\top}\right),\mbox{Prod}\left(\mathbf{Q},\mathbf{D},\mathbf{D}^{\top}\right)^{\top^{2}},\mbox{Prod}\left(\mathbf{Q},\mathbf{D},\mathbf{D}^{\top}\right)^{\top}\right),
\]
\[
\mbox{such that}
\]
\[
\mbox{Prod}\left(\mathbf{Q},\mathbf{Q}^{\top^{2}},\mathbf{Q}^{\top}\right)=\boldsymbol{\Delta},
\]
\[
\left[\mathbf{D}\right]_{ijk}=\begin{cases}
\begin{array}{cc}
\lambda_{jk}=\lambda_{kj}>0 & \mbox{if }0\le i=k<n\\
0 & \mbox{otherwise }
\end{array}\end{cases},
\]
where $\mathbf{Q}\in\mathbb{R}^{n\times n\times n}$ and $\mathbf{D}\in\mathbb{R}_{\ge0}^{n\times n\times n}$,
then 
\[
\forall\:0\le k<n,\quad\left(\mathbf{M}_{k}\left(\mathbf{z}\right)\right)^{\top}=\mathbf{M}_{k}\left(\mathbf{z}\right).
\]
Consequently each $\mathbf{M}_{k}\left(\mathbf{z}\right)$ is diagonalizable
and has real eigenvalues for all choices of the vector $\mathbf{z}$.
In particular, for $n=2$ it suffices to choose $\mathbf{z}$ such
that
\[
\forall\:0\le k<2,\quad\mbox{Tr}\left(\mathbf{M}_{k}\left(\mathbf{z}\right)\right)\ge0,\quad\det\left(\mathbf{M}_{k}\left(\mathbf{z}\right)\right)\ge0,
\]
which asserts that 
\[
q_{000}^{3}z_{0}+q_{001}q_{100}q_{101}z_{0}+q_{000}q_{001}q_{100}z_{1}+q_{101}^{3}z_{1}\ge0,
\]
\[
\left(q_{000}^{3}z_{0}+q_{000}q_{001}q_{100}z_{1}\right)\left(q_{001}q_{100}q_{101}z_{0}+q_{101}^{3}z_{1}\right)-\left(q_{000}q_{001}q_{100}z_{0}+q_{001}q_{100}q_{101}z_{1}\right)^{2}\ge0,
\]
\[
q_{010}^{3}z_{0}+q_{011}q_{110}q_{111}z_{0}+q_{010}q_{011}q_{110}z_{1}+q_{111}^{3}z_{1}\ge0,
\]
\[
\left(q_{010}^{3}z_{0}+q_{010}q_{011}q_{110}z_{1}\right)\left(q_{011}q_{110}q_{111}z_{0}+q_{111}^{3}z_{1}\right)-\left(q_{010}q_{011}q_{110}z_{0}+q_{011}q_{110}q_{111}z_{1}\right)^{2}\ge0.
\]
For all $\mathbf{x}$, $\mathbf{y}$ and $\mathbf{z}$ chosen as indicated
above we have 
\[
\min_{0\le i,j,k,t<n}\left(\lambda_{it}^{2}\,\lambda_{jt}^{2}\,\lambda_{kt}^{2}\right)\le\frac{\mbox{Prod}_{\mathbf{A}}\left(\mathbf{x}^{T^{2}},\mathbf{y}^{T},\mathbf{z}\right)}{\mbox{Prod}\left(\mathbf{x}^{T^{2}},\mathbf{y}^{T},\mathbf{z}\right)}\le\max_{0\le i,j,k,t<n}\left(\lambda_{it}^{2}\,\lambda_{jt}^{2}\,\lambda_{kt}^{2}\right).
\]

\section{Some related algorithmic problems}

\subsection{Logarithmic least square}

Let $\mathbf{A}\in\mathbb{C}^{m\times n}$, $\mathbf{b}\in\mathbb{C}^{m\times1}$
and consider the monomial constraints in the unknown $\mathbf{x}$
of size $n\times1$ vector 
\begin{equation}
\left\{ b_{i}=\prod_{0\le j<n}x_{j}^{a_{ij}}\right\} _{0\le i<m}.\label{Generic Monomial System}
\end{equation}
The \emph{logarithmic least square solution} to \eqref{Generic Monomial System}
is obtained by solving for $\mathbf{x}$ in the modified system 
\[
\left\{ \prod_{0\le t<m}b_{t}^{\overline{a_{ti}}}=\prod_{0\le j<n}\left(\prod_{0\le t<m}x_{j}^{\overline{a_{ti}}\,a_{tj}}\right)\right\} _{1\le i\le n}.
\]
By the least square argument the modified system is known to always
admit a solution vector $\mathbf{x}$ which minimizes 
\[
\sum_{0\le i<m}\left|\ln\left(b_{i}^{-1}\prod_{0\le j<n}x_{j}^{a_{ij}}\right)\right|^{2}.
\]
Such a solution is called the logarithmic least square solution of
the system \eqref{Generic Monomial System} and can be obtained via
the variant Gauss-Jordan elimination discussed in section 3.2.

\subsection{Logarithmic least square BM-rank one approximation}

Let $\rho$ denote some positive integer for which $0\le\rho\le n$.
A solution to the general BM-rank $\rho$ approximation of a cubic
$m$-th order hypermatrix $\mathbf{H}$ having side length $n$ is
obtained by solving for a BM conformable $m$-tuple $\left(\mathbf{X}^{(i)}\right)_{1\le i\le m}$
which minimize the norm

\[
\left\Vert \mathbf{H}-\sum_{0\le t<\rho}\mbox{Prod}_{\boldsymbol{\Delta}^{(t)}}\left(\mathbf{X}^{(1)},\mathbf{X}^{(2)},\cdots,\mathbf{X}^{(m)}\right)\right\Vert 
\]
Consequently, the BM-rank of $\mathbf{H}$ over $\mathbb{C}$ is the
smallest positive integer $\rho$ for which there exist a BM conformable
hypermatrix $m$-tuple $\left(\mathbf{X}^{(i)}\right)_{1\le i\le m}$
such that 
\[
\left\Vert \mathbf{H}-\sum_{0\le t<\rho}\mbox{Prod}_{\boldsymbol{\Delta}^{(t)}}\left(\mathbf{X}^{(1)},\,\mathbf{X}^{(2)},\,\cdots,\mathbf{X}^{(m)}\right)\right\Vert =0.
\]
It is easy to see that for all $m$-th order cubic hypermatrix $\mathbf{H}$
of side length $n$ 
\[
0\le\mbox{BM-rank}\left(\mathbf{H}\right)\le n.
\]

In particular constraints associated with the BM-rank $1$ problem
\[
\mathbf{H}=\mbox{Prod}_{\boldsymbol{\Delta}^{(0)}}\left(\mathbf{X}^{(1)},\mathbf{X}^{(2)},\cdots,\mathbf{X}^{(m)}\right),
\]
are monomial constraints of same type as the ones in \eqref{Generic Monomial System}.
The corresponding system admits no solution if BM-rank$\left(\mathbf{H}\right)>1$.
Our proposed BM-rank $1$ approximation of $\mathbf{H}$ is thus obtained
by solving the constraints in the logarithmic least square sense.

\subsection{Logarithmic least square direct sum and Kronecker product approximation}

Let $\mathbf{A}$ denote a cubic $m$-th order hypermatrix of side
length $n$ such that 
\[
\mathbf{A}=\bigoplus_{1\le j\le\beta}\mathbf{A}^{(j)},
\]
where $\mathbf{A}^{(j)}\in\mathbb{C}^{2^{j}\times2^{j}\times\cdots\times2^{j}}$.
A direct sum and Kronecker product approximation of $\mathbf{A}$
is obtained by solving for entries of a hypermatrix $\mathbf{B}$
subject to two constraints. The first constraints asserts that $\mathbf{B}$
must be generated by some a arbitrary combinations of Kronecker products
and direct sums of cubic side length $2$ hypermatrices. The second
constraint asserts that $\mathbf{B}$ should be chosen so as to minimize
the norm $\left\Vert \mathbf{A}-\mathbf{B}\right\Vert $. The problem
reduces to a system of the same form as \eqref{Generic Monomial System}
and is given by 
\[
\left\{ \mathbf{A}^{(j)}=\bigotimes_{0\le i<j}\mathbf{X}^{(i,j)}\right\} _{1\le j\le\beta},
\]
where $\mathbf{X}^{(i,j)}\in\mathbb{C}^{2\times2\times\cdots\times2}$.
Consequently the system admits no solution if $\mathbf{A}$ is not
generated by a combination of Kronecker product and direct sums of
side length $2$ hypermatrices. Our proposed direct sum and Kronecker
product approximation of $\mathbf{A}$ is obtained by solving the
corresponding system in the logarithmic least square sense.

\bibliographystyle{amsalpha}
\bibliography{mybib}

\end{document}